\documentclass[11pt,a4paper,reqno]{amsart}
\usepackage{amsfonts,amsthm, amscd, epsfig, amsmath, amssymb,enumerate}

\title[Random walks, exclusion processes  on  random infinite clusters]{Random walks and  exclusion processes among random
conductances on random infinite clusters: homogenization and
hydrodynamic limit}

\author{Alessandra Faggionato}
\address{Alessandra Faggionato. Dipartimento di Matematica ``G. Castelnuovo", Universit\`a   ``La
  Sapienza''. P.le Aldo Moro  2, 00185  Roma, Italy. e--mail:
  faggiona@mat.uniroma1.it}

\numberwithin{equation}{section}

\DeclareMathSymbol{\leqslant}{\mathalpha}{AMSa}{"36} 
\DeclareMathSymbol{\geqslant}{\mathalpha}{AMSa}{"3E} 
\DeclareMathSymbol{\eset}{\mathalpha}{AMSb}{"3F}     
\renewcommand{\leq}{\;\leqslant\;}                   
\renewcommand{\geq}{\;\geqslant\;}                   

\newcommand{\be}{\begin{equation}}

\def\1{\ifmmode {1\hskip -3pt \rm{I}} \else {\hbox {$1\hskip -3pt \rm{I}$}}\fi}

\newtheorem{Th}{Theorem}[section]
\newtheorem{Le}[Th]{Lemma}
\newtheorem{Pro}[Th]{Proposition}
\newtheorem{Cor}[Th]{Corollary}

\newtheorem{Def}{Definition}

\newcommand{\cA}{\ensuremath{\mathcal A}}
\newcommand{\cB}{\ensuremath{\mathcal B}}
\newcommand{\cC}{\ensuremath{\mathcal C}}
\newcommand{\cD}{\ensuremath{\mathcal D}}
\newcommand{\cE}{\ensuremath{\mathcal E}}
\newcommand{\cF}{\ensuremath{\mathcal F}}
\newcommand{\cG}{\ensuremath{\mathcal G}}

\newcommand{\cL}{\ensuremath{\mathcal L}}
\newcommand{\cM}{\ensuremath{\mathcal M}}

\newcommand{\cP}{\ensuremath{\mathcal P}}

\newcommand{\bbE}{{\ensuremath{\mathbb E}} }

\newcommand{\bbI}{{\ensuremath{\mathbb I}} }

\newcommand{\bbL}{{\ensuremath{\mathbb L}} }

\newcommand{\bbP}{{\ensuremath{\mathbb P}} }
\newcommand{\bbQ}{{\ensuremath{\mathbb Q}} }
\newcommand{\bbR}{{\ensuremath{\mathbb R}} }

\newcommand{\bbX}{{\ensuremath{\mathbb X}} }

\newcommand{\bbZ}{{\ensuremath{\mathbb Z}} }

\newcommand{\oo}{{\tilde \omega} }
\newcommand{\mmu}{{\bf \mu}}

 \newcommand{\strano}{\rightharpoonup}
    \let\d=\delta  \let\e=\varepsilon
        \let\l=\lambda
      \let\o=\omega    \let\p=\pi  
  \let\s=\sigma \let\t=\tau   
  \let\z=\zeta
   \let\G=\Gamma  \let\L=\Lambda 
\let\O=\Omega      

\def\\{\hfill\break}

\def\tthsp{\kern .083333 em}

\def\?{\mskip -10mu}
\def\indbox#1{\hbox to \parindent{\hfil\ #1\hfil} }

\newfam\msafam
\newfam\msbfam
\newfam\eufmfam
\def\hexnumber#1{%
  \ifcase#1 0\or 1\or 2\or 3\or 4\or 5\or 6\or 7\or 8\or
  9\or A\or B\or C\or D\or E\or F\fi}
\font\tenmsa=msam10 \font\sevenmsa=msam7 \font\fivemsa=msam5
\textfont\msafam=\tenmsa \scriptfont\msafam=\sevenmsa
\scriptscriptfont\msafam=\fivemsa
\edef\msafamhexnumber{\hexnumber\msafam}%
\mathchardef\restriction"1\msafamhexnumber16 \mathchardef\ssim"0218
\mathchardef\square"0\msafamhexnumber03
\mathchardef\eqd"3\msafamhexnumber2C
\def\QED{\ifhmode\unskip\nobreak\fi\quad
  \ifmmode\square\else$\square$\fi}
\font\tenmsb=msbm10 \font\sevenmsb=msbm7 \font\fivemsb=msbm5
\textfont\msbfam=\tenmsb \scriptfont\msbfam=\sevenmsb
\scriptscriptfont\msbfam=\fivemsb
\def\Bbb#1{\fam\msbfam\relax#1}
\font\teneufm=eufm10 \font\seveneufm=eufm7 \font\fiveeufm=eufm5
\textfont\eufmfam=\teneufm \scriptfont\eufmfam=\seveneufm
\scriptscriptfont\eufmfam=\fiveeufm

\def\({\left(}
\def\){\right)}
\let\Pp=\bP

\let\neper=e
\let\ii=i

\let\<=\langle
\let\>=\rangle

\def\PP{ \mathop\Pp\nolimits }

\outer\def\nproclaim#1 [#2]#3. #4\par{\medbreak \noindent
   \talato(#2){\bf #1 \Thm[#2]#3.\enspace }%
   {\sl #4\par }\ifdim \lastskip <\medskipamount
   \removelastskip \penalty 55\medskip \fi}
\def\thmm[#1]{#1}
\def\teo[#1]{#1}
\def\sttilde#1{%
\dimen2=\fontdimen5\textfont0 \setbox0=\hbox{$\mathchar"7E$}
\setbox1=\hbox{$\scriptstyle #1$} \dimen0=\wd0 \dimen1=\wd1
\advance\dimen1 by -\dimen0 \divide\dimen1 by 2
\vbox{\offinterlineskip%
   \moveright\dimen1 \box0 \kern - \dimen2\box1}
}
\begin{document}

\maketitle
\begin{abstract}
We consider a stationary and ergodic random field $\{\o (b):b \in
\bbE_d  \}$ parameterized by the family of  bonds in $\bbZ^d$,
$d\geq 2$. The random variable $\o(b)$  is  thought of as the
conductance of  bond $b$ and  it ranges in a finite interval
$[0,c_0]$. Assuming that the set of bonds with positive conductance
has a unique infinite cluster $\cC(\o)$, we prove homogenization
results for the  random walk  among random conductances on
$\cC(\o)$.
 As a byproduct, applying the general criterion of \cite{F} leading
to the hydrodynamic limit of exclusion processes with
bond--dependent transition  rates,  for almost all realizations of
the environment we prove  the hydrodynamic limit of simple exclusion
processes among  random conductances   on  $\cC(\o)$. The
hydrodynamic equation is given by a heat equation whose diffusion
matrix does not depend on the environment. We do not require any
ellipticity condition. As special case, $\cC(\o)$ can be the
infinite cluster of supercritical Bernoulli bond percolation.

\medskip

 \noindent \emph{Key words}: disordered system,
 bond percolation,  random walk in
random environment, exclusion process, homogenization.

\smallskip
\noindent \emph{AMS 2000 subject classification}:  60K35, 60J27,
82C44.
\end{abstract}

\section{Introduction}

We consider a stationary and ergodic random field
$\o=\bigl(\o(b)\,:\, b \in \bbE_d\bigr)$, parameterized by the set
$\bbE_d$ of non--oriented bonds in $\bbZ^d$, $d\geq 2$, such that
$\o(b) \in [0,c_0]$ for some fixed positive constant $c_0$. We call
 $\o$ the  conductance field and we interpret $\o(b)$ as the
conductance at bond $b$.  We assume  that the
 network of bonds $b$ with positive conductance  has a.s.  a unique infinite cluster
 $\cC$, and call $\cE$ the associated bonds. Finally, we consider
  the exclusion process on the graph
 $(\cC, \cE)$ with
 generator $\bbL $  defined on local functions $f$ as
 $$ \bbL f (\eta ) = \sum _{b \in \cE } \o(b) \bigl( f(\eta ^b ) -
 f(\eta) \bigr)\,, \qquad \eta \in \{0,1\}^{\cC   }
\,,  $$
 where  the
 configuration $\eta^b $ is obtained from  $\eta$ by exchanging the
 values of $\eta_x$ and  $\eta_y$, if $b = \{x,y\}$.
 If   $\bigl(\o(b) \,:\,b \in \cE \bigr) $ is  a family of i.i.d.  random
   variables such that  $\o(b)$   is positive with probability $p$ larger than the
   critical threshold $p_c$ for Bernoulli bond percolation,
    then the above process is an exclusion process  among positive random conductances on the
   supercritical percolation cluster.
 If  $p=1$, then   the  model reduces to
the exclusion process among positive random conductances on
$\bbZ^d$.

\smallskip
Due to the disorder,  the above model is an example of non--gradient
exclusion process, in the sense that the transition rates cannot be
written as gradient of some local function on $\{0,1\}^{\cC}$
\cite{KL} (with exception of the case of constant conductances).
Despite this fact, the hydrodynamic limit of the exclusion process
can be proven without using the very sophisticated techniques
developed for non--gradient systems  (cf. \cite{KL} and references
therein), which in addition would require non trivial spectral gap
estimates  that fail in the case of conductances non bounded from
below by a positive constant (cf. Section 1.5 in \cite{M}). The
strong simplification comes from the fact that, since the transition
rates depend only on the bonds but not on the particle
configuration, the function $\bbL \eta_x$, where $\eta_x$ is the
occupancy number at site $x\in \cC $, is a linear combinations of
occupancy numbers.
 Due to
 this degree conservation  the analysis of the limiting behavior
 of the random  empirical measure $\p(\eta)  = \sum _{x\in \cC}\eta_x  \d_x $  is
 strongly    simplified w.r.t.
 disordered models with transition rates depending both on the disorder
 and on the particle configuration \cite{Q1}, \cite{FM}, \cite{Q2}.
Moreover,  the function $\bbL \eta_x$ can be written as $(\cL
\eta)_x$,  where $\cL$ is the generator of the random walk on $\cC$
of  a single particle among the random conductances $\o(b)$ and
$\eta\in \{0,1\}^{\cC }$ is thought of  as an observable on the
state space $\cC $ of the random walk. This observation allows to
derive  the hydrodynamic limit of the exclusion process on $\cC$
from  homogenization results for the random walk on $\cC$. This
reduction has been performed  in \cite{N} for the exclusion process
on $\bbZ$ with bond--dependent conductances,  the  method has been
improved  and extended to the $d$--dimensional case in \cite{F}. The
arguments followed in \cite{F} are very general and can be applied
also to exclusion processes with bond--dependent transition rates on
\emph{general (non--oriented) graphs}, even with   \emph{non
diffusive behavior} (see \cite{FJL} for an example of application).
The reduction to a homogenization problem can be performed also by
means of the method of corrected empirical measure, developed in
\cite{JL} and \cite{J}. In \cite{JL} the authors consider exclusion
processes on $\bbZ$ with bond--dependent rates, while  in \cite{J}
the author proves the hydrodynamic limit for exclusion processes
with bond--dependent rates on triangulated domains and  on the
Sierpinski gasket. Moreover, in \cite{J} the author reobtains the
hydrodynamic limit of the exclusion process on $\bbZ^d$ among
 conductances bounded from above and from  below by  positive constants. Note that
this last result follows at once by applying the standard
non--gradient methods (in this case, their application becomes
trivial) or the discussion given in \cite{F}[Section 4]. Moreover,
it is reobtained in the present paper by taking  $(\o_b\,:\, b \in
\bbE_d)$  as independent strictly positive random variables.

\smallskip

By the above methods \cite{N}, \cite{F}, \cite{JL}, \cite{J}, the
proof of the hydrodynamic limit of exclusion processes on graphs
with bond--dependent rates  reduces to  a
 homogenization problem. In
our contest,
 we solve this problem by means of the notion of
\emph{two--scale convergence}, which is particularly fruitful when
dealing with homogenization problems on singular structures.
The notion of two--scale convergence was introduced by G. Nguetseng
\cite{Nu} and developed by G. Allaire \cite{A}.  In particular, our
proof is inspired by the method developed in \cite{ZP} for
differential operators on singular structures. Due to the ergodicity
and the $\bbZ^d$--translation invariance of the conductance field,
the arguments of \cite{ZP} can be simplified: for example,  as
already noted in \cite{MP}, one can avoid   the introduction of the
Palm distribution. Moreover, despite \cite{ZP} and previous results
of homogenization of random walks in random environment (see
\cite{Ko}, \cite{Ku}, \cite{PR} for example), in the present setting
we are able to avoid ellipticity assumptions.

\smallskip

We point out that recently the quenched central limit theorem for
the   random walk among constant conductances  on the supercritical
percolation cluster has been proven in \cite{BB} and \cite{MP}, and
previously for dimension $d \geq 4 $ in \cite{SS}. Afterwards, this
result has been extended to the case of i.i.d. positive bounded
conductances on the supercritical percolation cluster  in \cite{BP}
and \cite{M}. Their proofs  are very robust and use sophisticated
techniques and estimates (as heat kernel estimates, isoperimetric
estimates, non trivial percolation results, ...), previously
obtained in other papers. In the case of  i.i.d. positive bounded
conductances
 on the supercritical percolation cluster, we did not try to derive our homogenization results from the above
quenched CLTs (this route would anyway require some technical work).
Usually, the proof of the quenched CLT is based on homogenization
ideas and not viceversa. And in fact, our proof of homogenization is
much simpler than the above proofs of the quenched CLT. Hence, the
strategy we have followed has the advantage to be self--contained,
simple and to give a genuine homogenization result without
ellipticity assumptions. This would have not been achieved choosing
the above alternative route. In addition, our results cover more
general random conductance fields, possibly with correlations, for
which a proof of the quenched  CLT for the associated random walk is
still lacking.

\smallskip

The paper  is organized as follows: In Section \ref{natalino} we
give a more detailed description of the exclusion process and the
random walk among random conductances on the infinite cluster $\cC$.
In addition,  we state our main results concerning the hydrodynamic
behavior of the exclusion process (Theorem \ref{annibale})  and the
homogenization of the random walk (Theorem \ref{prop_hom} and
Corollary \ref{tg5}). In Section \ref{aleale} we show how to apply
the results of \cite{F} in order to derive Theorem \ref{annibale}
from Corollary  \ref{tg5}, while the remaining sections are focused
on the homogenization problem. In particular, the proof of Theorem
\ref{prop_hom} is given in Section \ref{kukakuka}, while the proof
of Corollary \ref{tg5} is given in Section \ref{kukakukabis}.
Finally, in the Appendix we prove Lemma \ref{silurino}, assuring
that the class of random conductance fields satisfying our technical
assumptions is large.

\section{Models and results}\label{natalino}

\subsection{The environment}\label{salutare}
The   environment  modeling the disordered medium is given by a
stationary and ergodic random field $\o=\bigl(\o(b)\,:\, b \in
\bbE_d\bigr)$, parameterized by the set $\bbE_d$ of non--oriented
bonds in $\bbZ^d$, $d\geq 2$. Stationarity and ergodicity refer to
the natural  action of the group of $\bbZ^d$--translations.
 $\o$ and $\o(b)$ are  thought of as the  conductance field  and the
conductance at bond $b$, respectively.  We call $\bbQ$ the law of
the field $\o$ and we assume that

\bigskip
 (H1) \centerline{$ \o(b)
\in [0,c_0]\,, \qquad \text{$\bbQ$--a.s.}    \qquad \text{$  $}$}

\bigskip

\noindent for some fixed positive constant $c_0$. Hence, without
loss of generality, we can suppose that $\bbQ$ is a  probability
measure on the product space $\Omega:=[0,c_0] ^{\bbE _d} $.
Moreover, in order to simplify the notation, we write $\o(x,y)$ for
the conductance $\o(b)$ if $b=\{x,y\}$. Note that $\o(x,y)=\o(y,x)$.

\smallskip

  Consider the random graph $G (\o)= \bigl(  V (\o), E
(\o) \bigr)$ with  vertex set $V (\o)$ and bond set $E(\o)$ defined
as
\begin{align*}
& E (\o):=\bigl\{b\in \bbE_d\,:\, \o (b) >0\bigr \}\,, \\
& V (\o) :=\bigl\{x\in \bbZ ^d \,:\, x\in b \text{ for some } b \in
 E (\o ) \bigr\}.
\end{align*}
Due to ergodicity,  the translation invariant Borel subset
$\O_0\subset \O$ given by the configurations $\o$ for which the
graph $ G (\o)$ has a unique infinite connected component (cluster)
$\cC (\o) \subset  V (\o )$  has $\bbQ$--probability $0$ or $1$. We
assume that

\bigskip

  (H2)\centerline{$\bbQ(\O_0)=1.    \qquad \text{$  $}$}

\bigskip

\noindent
 Below, we denote by $\cE (\o)$ the bonds in $E(\o)$
connecting points of $\cC(\o)$ and we will often understand the fact
that $\o \in \O_0$.


\smallskip

Define $B (\O)$ as the family of bounded Borel functions on $\O$ and
let  $\cD$ be the $d\times d $ symmetric matrix characterized by the
variational formula
\begin{equation}\label{varcar}
(a, \cD a) =\frac{1}{m} \inf _{\psi \in B (\O) }\left\{
 \sum _{e\in \cB_*} \int _\O \o (0,e)  ( a_e+ \psi (\t_e\o)-\psi (\o) ) ^2 \bbI _{0,e\in \cC (\o) }
  \bbQ ( d\o)\right\} \;,\;\;\forall a\in \bbR^d \,,
  \end{equation}
where $\cB_*$ denotes the canonical basis of $\bbZ^d$,
\begin{equation}\label{defm} m:= \bbQ \left(0 \in \cC (\o)
\right)\,  \end{equation} and the translated environment  $\t_e \o $
is defined as $\t_e \o (x,y)=\o(x+e, y+e)$ for all bonds $\{x,y\}$
in $\bbE_d$. In general, $\bbI_A$ denotes the characteristic
function of $A$. Our last assumption on $\bbQ$ is that the matrix
$\cD$ is strictly positive:

\medskip

\medskip
(H3) \centerline{ $(a, \cD a)>0\,, \qquad \forall a \in \bbR^d\,:\,
a \not =0 \,. \qquad \text{$  $}$}

\medskip

\medskip

\noindent In conclusion, our hypotheses on the random field $\o$ are
given by stationarity, ergodicity, (H1),(H2) and (H3). The lemma
below shows that they are fulfilled by a large class of random
fields. In order to state it, given $c>0$  we define the random
field $\hat \o_c=\bigl(\hat \o_c(b)\,:\, b \in \bbE_d\bigr)$  as
 \begin{equation}\label{paranza3}
\hat \o_c (b) = \begin{cases} 1 & \text{ if } \o(b)>c\,, \\
0 & \text{ otherwise}\,.
\end{cases}
\end{equation}
For $c=0$ we simply set $\hat \o:=\hat \o _0$.

\begin{Le}\label{silurino}
Hypotheses (H2) and (H3) are satisfied  if there exists a positive
constant $c$ such that the random field $\hat \o _c$ stochastically
dominates a supercritical Bernoulli bond percolation. In particular,
if $\hat\o  $ itself is a supercritical Bernoulli bond percolation,
then (H2) and (H3) are verified.

If  the field $\o$ is reflection invariant or isotropic (invariant
w.r.t. $\bbZ^d$ rotations by $\p/2$), then $\cD$ is a diagonal
matrix. In the isotropic case $\cD$ is a multiple of the identity.
In particular, if $\o$ is given by i.i.d. random conductances
$\o(b)$, then $\cD$ is a  multiple of the identity.
\end{Le}

We postpone the proof of the above Lemma to the Appendix.

\subsection{The exclusion process on the infinite
cluster $\cC(\o)$}  Given a realization $\o$ of the environment, we
consider the  exclusion process  $\eta (t) $ on the graph  $\cG
(\o)= \bigl( \cC (\o), \cE (\o)\bigr) $ with   exchange rate $\o(b)$
at bond $b$. This is the
 Markov process with paths $\eta(t)$ in the Skohorod space
$D\bigl( [0,\infty), \{0,1\} ^{\cC (\o) } \bigr)$ (cf. \cite{B})
whose Markov generator $\bbL _\o $ acts on local functions as
\begin{equation}
\bbL_\o f (\eta) = \sum _{e\in \cB_*} \sum _{\substack{x\in
\cC(\o)\,:\,\\
x+e \in \cC (\o) }  } \o(x,x+e)
 \left( f(\eta ^{x,x+e})- f(\eta)\right)\,,
\end{equation}
where  in general
$$ \eta^{x,y} _z  =
\begin{cases}
\eta_y,&\text{ if }  z=x\,,\\
\eta_x, &\text{ if }z=y\,,\\
\eta_z, &\text{ if } z \neq x,y\,.\\
\end{cases}
$$
We recall that a function $f$ is called local if $f(\eta)$  depends
only on $\eta_x$  for a finite number of sites $x$. By   standard
methods \cite{L} one can prove that the above exclusion process
$\eta(t)$ is well defined.

\smallskip

 Every configuration $\eta $ in the state space $\{0,1\}^{\cC (\o)}$
 corresponds to a system of particles on $\cC (\o) $ if one considers a
site $x$ occupied by a particle if $\eta_x=1$ and vacant if $\eta_x
=0$. Then the exclusion process is given by a stochastic dynamics
where particles can lie only on sites $x\in \cC (\o)$ and can jump
from the original site $x$ to the vacant site $y\in \cC (\o)$ only
if the bond $\{x,y\}$ has positive conductance, i.e. $x$ and  $y$
are connected by a bond in $\cG (\o)$. Roughly speaking, the
dynamics can be described as follows: To each bond $b=\{x,y\} \in
\cE (\o)$ associate an exponential alarm clock with mean waiting
time  $1/\o (b)$. When the clock rings,  the particle configurations
at sites $x$ and $y$ are exchanged and   the alarm clock restarts
afresh. By Harris'  percolation argument \cite{D}, this construction
can be suitably formalized. Finally, we point out that the only
interaction between particles is given by site exclusion.

\smallskip



We can finally describe the hydrodynamic limit of the above
exclusion process  among random conductances $\o(b)$ on the infinite
 cluster $\cC (\o)$.  If the initial distribution is given by the
probability measure $\mu$ on $\{0,1\}^{\cC (\o)}$,  we denote by
$\bbP _{\o, \mu}$ the law of the resulting  exclusion process.

\begin{Th}\label{annibale}
For $\bbQ$  almost all environments $\o$ the following holds. Let
$\rho _0: \bbR ^d \rightarrow [0,1]$ be a Borel function and let
$\{\mu _\e\} _{\e>0}$ be a family of probability measures on
$\{0,1\} ^{\cC (\o)} $ such that, for all $\d>0$ and all real
functions $\varphi$ on $ \bbR^d $   with  compact support (shortly
$\varphi  \in C_c (\bbR ^d)) $, it holds
\begin{equation}\label{pasqua1}
\lim _{\e \downarrow 0 } \mu_\e \Big( \Big| \e^d  \sum _{x\in \cC
(\o) } \varphi (\e x )\,  \eta _x - \int _{\bbR ^d} \varphi (x)
\rho_0  (x) dx \Big|>\d \Big )=0\,.
\end{equation}
 Then, for all $t>0$, $\varphi \in C_c (\bbR ^d)$ and $\d>0$,
\begin{equation}\label{pasqua2}
\lim _{\e \downarrow 0 } \bbP _{\o, \mu_\e}  \Big ( \Big| \e^d \sum
_{x\in \cC (\o) } \varphi (\e x )\,  \eta _x (\e^{-2} t)  - \int
_{\bbR ^d} \varphi (x) \rho (x,t ) dx  \Big|>\d \Big )=0\,,
\end{equation}
where  $\rho: \bbR^d \times [0,\infty) \rightarrow \bbR$ solves the
heat equation
\begin{equation}\label{pasqua3}
\partial  _t \rho = \nabla \cdot ( \cD \nabla  \rho)=\sum_{i,j=1}^d
\cD _{i,j}  \partial^2_{x_i,x_j} \rho
\end{equation}
with boundary condition $\rho_0$ at $t=0$  and where the symmetric
matrix $\cD$ is variationally characterized by  (\ref{varcar}).
\end{Th}

 If a density profile $\rho_0$
can be approximated by a family of probability measures $\mu_\e$ on
$\{0,1\}^{\cC (\o) }$ (in the sense that  (\ref{pasqua1}) holds for
each $\d>0$ and $\varphi\in C_c (\bbR^d)$),
  then it must be  $0 \leq \rho_0 \leq m $ a.s.
On the other hand, if $\rho_0:\bbR ^d \rightarrow [0,m] $ is a
Riemann integrable function, then it is simple to exhibit for
$\bbQ$--a.a. $\o$  a family of probability measures $\mu_\e$ on
$\{0,1\}^{\cC (\o) }$  approximating $\rho_0$. To this aim we
observe that,
 due to the ergodicity of $\bbQ$ and by separability arguments,
    for $\bbQ$ a.a. $\o$ it holds
\begin{equation}\label{scudetto}
\lim_{\e\downarrow 0} \e^d  \sum _{x\in \cC (\o) } \varphi (\e x )=
m \int_{\bbR ^d} \varphi (x) dx \,,
\end{equation}
 for each Riemann integrable  function $\varphi : \bbR^d
\rightarrow \bbR$ with compact support. Fix such an environment
$\o$. Then, it is enough to define $\mu_\e$ as the unique product
probability measure on $\{0,1\}^{\cC (\o ) }$ such that $\mu_\e
(\eta_x=1) = \rho_0 (\e x )/m $ for each $x \in \cC (\o)$. Since the
random variable $\e^d \sum _{x\in \cC (\o) } \varphi (\e x )\, \eta
_x $ is the  sum of independent random variables,  it is simple to
verify that its mean equals $\e^d \sum _{x\in \cC (\o) } \varphi (\e
x )\rho_0 (\e x)/m $ and its  variance equals  $\e^{2d} \sum _{x\in
\cC (\o) } \varphi^2 (\e x )[\rho_0 (\e x)/m][1-\rho_0(\e x) /m] $.
The thesis then follows by means of (\ref{scudetto}) and the
Chebyshev inequality.
\medskip

 The proof of Theorem \ref{annibale} is  given in Section
\ref{aleale}. As already mentioned, it is  based on the general
criterion for the hydrodynamic limit of exclusion processes with
bond--dependent transition  rates, obtained  in \cite{F} by
generalizing an argument of \cite{N}, and homogenization results for
the random walk on $\cC (\o)$ with jump rates $\o(b)$, $b \in \cE
(\o)$, described below.

\subsection{The random walk among random conductances on the
infinite cluster $\cC (\o) $}

Given $\o \in \O$ we denote by $X_\o(t|x)$  the continuous--time
random walk on $\cC (\omega)$ starting at $x\in \cC (\o) $, whose
Markov generator $\cL _\omega $ acts  on bounded functions
$g:\cC(\o) \rightarrow \bbR$ as
\begin{equation}
 \cL _\omega g (x) = \sum _{\substack{ y\,:\, y\in \cC(\omega) \\
|x-y|=1} } \o (x,y) \left(g(y)- g(x)\right)\,, \qquad x\in \cC
(\omega)\,.
 \end{equation}
The dynamics can be described
 as follows. After arriving at site $z\in \cC (\o) $, the particle
  waits an exponential time of parameter
$$ \l_\o (z):= \sum  _{\substack{ y\,:\, y\in \cC(\omega) \\
|z-y|=1} } \o (z,y)
$$
and then jumps to a site $y \in \cC (\o)$,  $|z-y|=1$, with
probability $\o (z,y) / \l_\o  (z)$. Since the jump rates are
symmetric, the counting measure on $\cC (\o)$ is reversible for the
random walk.

In what follows, given $\e>0$ we will consider the rescaled random
walk \begin{equation}\label{classica} X_{\e , \o}(t|x) = \e X _\o (
\e ^{-2} t | \e ^{-1} x) \end{equation}
 with starting point $x\in \e
\cC (\o )$. We denote by $\mu _\o ^\e$ the reversible  rescaled
counting measure
  $$
  \mu _\omega ^\e =\e^d \sum _{x\in \cC (\o) } \delta _{\e x}
  $$
and write  $ \cL ^\e _\o$ for the symmetric operator   on $L^2 (\mu
^\e _\o) $  defined   as
\begin{equation}
 \cL ^\e_\omega g (\e x) = \e ^{-2}\sum _{\substack{ y\,:\, y\in \cC(\omega) \\
|x-y|=1} }\o (x,y )  \left(g(\e y)- g(\e x)\right)\,, \qquad x\in
\cC (\omega)\,.
 \end{equation}

\smallskip

Due to (\ref{scudetto}),   for almost all $\o \in \O$ the  measure
$\mu^\e _\o $ converges vaguely to the measure $m \, dx $, where the
positive constant $m$ is defined in (\ref{defm}).   In what follows,
$\|\cdot \|_{\mu_\o^\e}$ and $ ( \cdot, \cdot)_{\mu_\o^\e}$ will
denote the norm and the inner product in $L^2 (\mu_\o^\e)$,
respectively.
 We recall a standard definition in homogenization theory (cf. \cite{Z}, \cite{ZP} and reference
 therein):
\begin{Def} Fix $\o \in \O_0$. Given a family  of functions $f^\e_\o  \in L^2(\mu^\e_\o)$ parameterized
by $\e>0$
 and a
function $f \in L^2 ( m\,dx )$,  $f^\e_\o $  weakly converges to $f$
(shortly, $f^\e _\o   \rightharpoonup f$)   if
\begin{equation}\label{rimmel}
\sup _\e \| f^\e_\o  \|_{\mu^\e _\o } <\infty,
\end{equation}
and
\begin{equation}
\lim _{\e\downarrow 0 } \int _{\bbR^d} f^\e _\o (x) \varphi (x)
\mu_\o^\e (dx) = \int _{\bbR^d} f (x) \varphi (x) m \, dx
\end{equation}
for all   functions $\varphi \in C_c ^\infty (\bbR^d)$, while
$f^\e_\o $ strongly converges to $f$ (shortly, $f^\e_\o \rightarrow
f$) if (\ref{rimmel}) holds and if
\begin{equation}\label{rimmelbis}
\lim _{\e\downarrow 0 } \int _{\bbR^d} f^\e_\o  (x) \varphi ^\e (x)
\mu_\o^\e (dx) =  \int _{\bbR^d} f  (x) \varphi (x) m \, dx\,,
\end{equation}
for every family $\varphi ^\e\in L^2 (\mu ^\e _\o)$ weakly
converging to $\varphi \in L^2 (m\, dx)$.
\end{Def}

 The strong convergence $f^\e_\o \rightarrow f$ admits the following
 characterization
  (cf. \cite{Z}[Proposition 1.1] and references therein):

 \begin{Le}\label{caldissimo}
Fix  $\o \in \O$ and functions $f^\e_\o  \in L^2 (\mu ^\e _\o )$, $f
\in L^2 (m dx)$, where $\e>0$. Then the strong convergence $ f^\e
_\o \rightarrow f$ is equivalent to  the weak convergence
 $f^\e _\o   \rightharpoonup f$  plus  the relation
 \begin{equation}
\lim_{\e \downarrow 0 } \int_{\bbR ^d} f^\e _\o (x) ^2 \mu ^\e _\o
(dx) =  m \int_{\bbR^d} f(x) ^2 dx \,.
 \end{equation}
 \end{Le}




We need now to isolate a Borel subset $\O_*\subset \O$ of {\sl
regular environments}. To this aim we first define $\O_1$ as the set
of $\o \in \O_0 $ (recall the definition of $\O_0$ given before
(H2)) such that
\begin{align}
& \lim _{\e\downarrow 0 } \mu^\e_\o (\L_\ell ) = m (2\ell)^d
\,,\qquad \L_\ell := [-\ell, \ell ]^d \,, \label{carlo}\\
& \lim _{\e\downarrow 0 } \int _{\bbR^d} \varphi (z) u
\bigl(\t_{z/\e} \o\bigr) \mu _\omega ^\e (dz) = \int _{\bbR^d}
\varphi (z) dz \int _\Omega  u(\o')\mmu (d\o') \,,\label{presilla}
\end{align}
for all $\ell>0$,  $ \varphi \in C_c(\bbR^d)$,  $u \in C (\O)$.
From the the ergodicity of $\bbQ$ and the separability of $C_c
(\bbR^d)$ and  $C ( \O)$, it is simple to derive that $\bbQ
(\O_1)=1$. The set of regular environments $\O_*$ will be defined in
Section \ref{picche}, after Lemma \ref{zorro}, since its definition
requires the concept of solenoidal forms.  We only mention here that
$\O_*\subset \O_1$ and $\bbQ(\O_*)=1$.

\smallskip

 We can finally state our main homogenization result,
similar to \cite{ZP}[Theorem 6.1]:
\begin{Th}\label{prop_hom} Fix $\o \in \O_*$. Let $f^\e_\o$ be a family of
 functions with $f^\e_\o\in L^2 (\mu^\e_\o)$ and let $f\in L^2 (m dx)$.
Given $\l>0$,  define $u^\e_\o\in L^2 (\mu^\e_\o)$, $u^0\in L^2 (m
dx)$ as the solutions of the following equations in $L^2
(\mu^\e_\o)$, $ L^2 (m dx)$  respectively:
\begin{align}
&\l u ^\e_\o -\cL ^\e _\o u^\e_\o = f^\e_\o  \,,\label{mammina}\\
& \l u^0- \nabla\cdot ( \cD \nabla   u^0) = f \,,\label{caldo}
\end{align}
where the symmetric matrix $\cD$ is variationally characterized in
(\ref{varcar}).

\smallskip

(i)  If
$f^\e _\o \strano  f $, then it holds
\begin{equation}\label{tommy}
 L^2 (\mu^\e _\o ) \ni u^\e_\o \strano  u^0 \in L^2 (m dx ) \,, \qquad \forall
 \l >0  \,.\end{equation}

\smallskip

(ii) If
$f^\e _\o \rightarrow f$, then it holds
\begin{equation}\label{tommy2}
 L^2 (\mu^\e _\o ) \ni u^\e_\o \rightarrow  u^0 \in L^2 (m dx )  \,, \qquad \forall
 \l >0     \,.\end{equation}

\smallskip

(iii) For each test function  $f \in C_c  (\bbR^d)$ and  setting
$f^\e_\o:=f $, it holds
\begin{equation}\label{tommybis}
\lim_{\e\downarrow 0} \int _{\bbR^d} \left| u^\e_\o (x) - u^0 (x)
\right| ^2 \mu_\o ^\e (dx) =0\, ,\qquad \forall \l>0\,.
\end{equation}
\end{Th}

The proof of Theorem \ref{prop_hom} will be given in Section
\ref{kukakuka}. We state  here an important  corollary  of the above
result:  Set  $P_{t,\o}^\e =e^{t \cL^\e_\o }$, $P_t =e ^{ t\nabla
\cdot (\cD\nabla\cdot) }$. Note that  $\bigl( P_{t,\o}^\e : t\geq 0
\bigr) $ is the $L^2 (\mu ^\e _\o )$--Markov semigroup associated to
the random walk $X_{\e,\o} (t|x)$, i.e.
\begin{equation}\label{semigruppo}
P_{t,\o}^\e g (x) = \bbE \left[ g \left(X_{\e,\o} (t|x)\right)
\right]\,, \qquad x \in \e \cC (\o )\,, g \in L^2 (\mu ^\e _\o )\,,
\end{equation}
while $P_t $ is the $L^2 (m dx  )$--Markov semigroup associated to
the diffusion with generator $\nabla \cdot (\cD \nabla \cdot)$.

As proven in Section \ref{kukakukabis}  it holds:
\begin{Cor}\label{tg5}
For each
$\o \in\O_*$,  given any function $f\in C _c (\bbR^d)$, it holds
\begin{equation}\label{kerbala1}
\lim _{\e\downarrow 0} \int _{\bbR^d} \bigl |P_{t,\o}^\e f (x) - P_t
f(x) \bigr|^2 \mu^\e_\o (dx) =0\,.
\end{equation}
In particular, for  each $\o \in \O_*$, given  any function $f\in C
_c (\bbR^d)$, it holds
\begin{equation}\label{kerbala2}
\lim _{\e\downarrow 0} \int _{\bbR^d} \bigl |P_{t,\o}^\e f (x) - P_t
f(x) \bigr| \mu^\e_\o (dx) =0\,.
\end{equation}
\end{Cor}


\section{Proof of Theorem \ref{annibale}}\label{aleale}
 As already mentioned, having the homogenization result given by
 Corollary \ref{tg5}, Theorem \ref{annibale} follows easily from
 the criterion of \cite{F}
 for the hydrodynamic limit of exclusion
 processes with bond--dependent  rates. The method discussed in \cite{F} is
 an improvement of the one developed in \cite{N} for the analysis of
bulk diffusion of 1d exclusion processes with bond--dependent rates.
Although in \cite{F} we have discussed the criterion with reference
to exclusion processes on $\bbZ^d$, as the reader can check the
method  is very general and can be applied to exclusion processes on
general graphs with bond--dependent rates, also under non diffusive
space--time rescaling and also when the   hydrodynamic behavior is
not    described by  heat equations (cf. \cite{FJL} for an example).

\smallskip

 The following proposition is the main technical tool in
order to reduce the proof  of the hydrodynamic limit to a problem of
homogenization for the random walk performed by a single particle
(in absence of other particles). Recall the definition
(\ref{semigruppo})  of the semigroup $P^\e_{t,\o}$ associated to the
rescaled random walk $X_{\e,\o}$ defined in (\ref{classica}).

\begin{Pro}\label{bolero} For $\bbQ$--a.a.  $\o$ the following holds.
Fix  $\d,t>0$, $\varphi \in C_c (\bbR^d)$ and let $\mu_\e$ be a
family  of probability measures on $\{0,1\}^{\cC (\o) }$. Then
\begin{equation}\label{vivace}
\lim_{\e\downarrow 0} \bbP _{\o, \mu_\e} \left(\,  \Big | \e^d
\sum_{x\in \cC (\o)} \varphi (\e x) \eta_x ( \e^{-2} t ) -\e^d
\sum_{x\in \cC (\o)  } \eta_x  (0) P^\e _{t, \o } \varphi (\e x )
\Big|>\d \, \right) =0\,.
\end{equation}
\end{Pro}

\begin{proof} One can prove the above proposition by the same arguments used
in \cite{F}[Section 3] or one can  directly  invoke the discussion
of \cite{F}[Section 4]  referred to exclusion processes on $\bbZ^d$
with non negative  transition rates, bounded from above. In fact, to
the probability measure  $\mu_\e$  on $\{0,1\}^{\cC (\o) }$  one can
associate   the  probability measure $\nu_\e $ on $\{0,1\}
^{\bbZ^d}$ so characterized: $\nu_\e$ is concentrated on the event
$$
\cA_\o  := \left \{ \eta \in \{0,1\}^{\bbZ^d} \,:\, \eta _x =0
\text{ if } x \not \in \cC (\o) \right\}\,,$$ and
$$
\nu_\e  \bigl(  \eta _x=1 \;\; \forall x \in \L  \bigr) = \mu_\e
\bigl(  \eta _x=1 \;\; \forall x \in \L  \bigr) \,, \qquad \forall
\L \subset \cC (\o) \,.
$$
Given  $\eta(t)\in \{0,1\}^{\cC(\o) }$,   define $ \s (t)\in
\{0,1\}^{\bbZ^d}$ as
$$
\s _x (t) =
\begin{cases} \eta_x (t) &\text{ if } x \in \cC (\o) \,;\\
0 & \text{ otherwise}\,.
\end{cases}
$$
Note that,  if $\eta(t)$ has law  $\PP _{\mu_\e}$, then  $\s(t)$ is
the exclusion process on $\bbZ^d$ with initial distribution $\nu_\e
$ and generator $f\rightarrow \sum _{b \in \bbE _d }\o (b ) \bigl (
f(\s^b)-f(\s)\bigr)$. In particular, Proposition \ref{bolero}
coincides with the limit (B.2) in \cite{F}[Section 4].
\end{proof}

We can now complete the proof of Theorem \ref{annibale}. First we
observe that whenever (\ref{pasqua1}) is satisfied for functions of
compact support, then it is satisfied also for functions that vanish
fast at infinity. Indeed, for our purposes it is enough to show that
the limit
\begin{equation}\label{tesina}
\lim_{\e \downarrow 0 }
 \mu_\e \Big(  \Big |
 \e^d \sum _{x\in \cC (\o) }
 f (\e x) \eta _x-\int_{\bbR^d}
 f  (x) \rho_0(x) dx \Big| >\d \Big)=0
\end{equation}
 is valid for all  functions $f  \in C(\bbR^d)$ such that
\begin{equation} \label{renania}
|f  (x)| \leq \frac{c}{1+|x|^{d+1} }\, , \qquad \forall x \in \bbR
^d \,.
\end{equation}
For such a function $f$, given  $\ell>0$ we can find $g_\ell \in C_c
(\bbR ^d)$ such that $g_\ell (x) = f (x) $ for all $x \in \bbR^d$
with $|x|\leq \ell$  and
 $| g_\ell (x)| \leq \frac{c}{1+|x|^{d+1}},$ for all $ x \in
 \bbR^d $.
 Then
 \begin{align}
 & \Big | \e^d \sum _{x\in \cC (\o) }f (\e x) \eta _x -
 \e^d \sum _{x\in \cC (\o) }
 g_\ell (\e x) \eta _x \Big | \leq \e^d \sum _{x\in \bbZ^d  \,:\,
 |\e x|>\ell
 }\frac{2c}{1+|\e x|^{d+1}} \leq c(\ell)\,,\label{flauto1}\\
 & \left| \int_{\bbR ^d} f  (x) \rho_0 (x) dx - \int_{\bbR^d}
 g_\ell (x) \rho_0 (x) dx \right| \leq  \int_{\{x\in \bbR^d\,:\,
 |x|>\ell\}}\frac{2c}{1+|x|^{d+1}}dx\leq c(\ell)\,,\label{flauto2}
\end{align}
for a suitable positive constant $c(\ell)$ going to zero as $\ell
\rightarrow \infty$. Since $g_\ell \in C_c (\bbR ^d)$,  by
assumption (\ref{pasqua1}) we obtain that
\begin{equation}
\lim_{\e \downarrow 0 }
 \mu_\e \Big(  \Big |
 \e^d \sum _{x\in \cC (\o) }
 g_\ell (\e x) \eta _x-\int_{\bbR^d}
 g_\ell (x)\rho_0(x)  dx \Big| >\d \Big)=0\,.
\end{equation}
The above limit together with (\ref{flauto1}) and (\ref{flauto2})
implies (\ref{tesina}) for all functions $f \in C(\bbR^d)$
satisfying (\ref{renania}).

In particular, (\ref{renania}) is valid for $f= P_t \varphi $, where
$P_t = e^{ t \nabla \cdot (\cD \nabla \cdot )  }$ and $\varphi \in C
_c (\bbR^d)$. Indeed, in this case $f$ decays exponentially. Due to
this observation, Proposition \ref{bolero} and the fact that
$$ \int _{\bbR ^d} P_t \varphi (x) \rho_0 (x) dx = \int_{\bbR^d}
\varphi (x) \, P_t \rho_0 (x) dx = \int_{\bbR^d} \varphi (x)\rho
(x,t) dx \,,$$
 in order
to prove (\ref{pasqua2}) it is enough to show that for $\bbQ $--a.a.
$\o$ it holds
\begin{equation}\label{salvezza}
\lim _{\e \downarrow 0 } \mu _\e \Big ( \e ^d \Big| \sum _{x \in \cC
(\o) } \eta_x P ^{\e} _{t,\o} \varphi  (\e x) - \sum _{x\in \cC (\o)
} \eta _x P _t \varphi  (\e x) \Big| > \d \Big) =0\,
\end{equation}
for any $\varphi  \in C  _c (\bbR^d)$ and $\d>0$.
 Since $$ \e ^d \Big| \sum _{x \in \cC (\o) } \eta_x P ^{\e}
_{t,\o} \varphi  (\e x) - \sum _{x\in \cC (\o) } \eta _x P _t
\varphi (\e x) \Big|\leq \int _{\bbR ^d} \bigl|  P ^{\e} _{t,\o}
\varphi ( x) - P _t \varphi  ( x)\bigr|\mu^\e _\o (dx)\,,
$$
(\ref{salvezza}) follows from (\ref{kerbala2}) of Corollary
\ref{tg5}. This concludes the proof of Theorem \ref{annibale}.

\section{Square integrable forms}\label{picche}

We now focus our attention on the proof of homogenization for the
random walk on the infinite cluster. To this aim, in this section we
introduce the Hilbert space of square integrable forms and show how
the variational formula (\ref{varcar}) can be interpreted in terms
of suitable orthogonal projections inside this Hilbert space.

\smallskip

 Let $\cM (\bbR^d)$ be the family of Borel measures on $\bbR^d$.
 Given $x\in \bbZ^d$, $y \in \bbR ^d$,  $\o\in \Omega$ and $\nu \in \cM
 (\bbR^d)$, $\t_x\o\in \O$ and $\t_y \nu \in \cM (\bbR^d)$ are
 defined as
 $$
\t_x \o (b) = \o (b+x) \;\; \forall b \in \bbE_d  , \qquad \t_y\nu
(A) =\nu  (A+y) \;\; \forall A\subset \bbR^d, \; \text{Borel set}\,.
$$
Note that
the family of random measures $\mu _\o^\e$ satisfies the identity
  \begin{equation}\label{fabri}
  \t_{\e x }\mu _\o ^\e = \mu _{\t_x \o} ^\e\,, \qquad \forall x \in
  \bbZ^d\,.
  \end{equation}
Let $\mmu$ be the measure  on $\O$  absolutely continuous w.r.t.
$\bbQ$ such that
\begin{equation}
\mmu (d\o )=
 \bbI _{0\in \cC (\o)
}\bbQ (d\o)
\end{equation}
 and define $\cB$ as
$$\cB:=\{e\in \bbZ^d\,:\, |e|=1 \}= \cB_*\cup (-\cB_*)
$$
($\cB_*$ is the set of coordinate vectors in $\bbR ^d$).

\smallskip

Given real functions $u $ defined on $  \O $ and  $v $ defined on
$\O \times \cB  $,   we  define   the {\sl gradient}
 $\nabla ^{(\o)} u: \O \times \cB \rightarrow \bbR$ and the {\sl  divergence}  $\nabla^{(\o)*} v: \O \rightarrow \bbR$, respectively, as
 follows:
\begin{align}
&  \nabla ^{(\o)} u (\o, e )=\hat \o (0,e)  \left[u (\t_e \o ) - u
(\o
)\right]\,, \qquad u :  \O \rightarrow \bbR\,, \\
& \nabla^{(\o)*} v (\o)= \sum _{e\in \cB} \o (0,e ) \left[
v(\o,e)-v(\t_e \o,-e)\right]\,, \qquad v: \O \times \cB \rightarrow
\bbR\,.
\end{align}
Moreover,  we endow the space $\O\times \cB$ with the Borel measure
$M$ defined by
$$
\int_{\O\times \cB} v dM= \sum _{e\in \cB} \int _\O  \o (0,e) v(\o,
e)\mmu (d\o)\,, $$ where $v$ is any  bounded Borel function on $\O
\times \cB$. Note that if $u \in L^p (\mu)$ and $v \in L^p (M)$ then
$\nabla ^{(\o) } u \in L^p (M)$ and $ \nabla^{(\o)*} v\in L^p
(\mu)$.

\smallskip

 The space $L^2(M)$ is called the space of {\sl square integrable forms}.
Note that $M$ gives zero measure to the set
$$
\bigl\{ (\o,e) \in \O \times \cB \,:\,
\{0,e\} \not \in \cE(\o)
 \bigr\}\,.
$$
Hence, given a square integrable form $v \in L^2 (M)$, we can always
assume that $v(\o,e)=0$ whenever $\{0,e\} \not \in \cE(\o) $.
  We define the space of
{\sl potential forms} $L^2_{\text{pot}}(M)$ and the space of {\sl
solenoidal forms} $L^2_{\text{sol}}(M) $ as follows:

\begin{Def} The space $L^2_\text{pot}(M)$ is the closure in $L^2 (M)$ of
the set of gradients $\nabla ^{(\o)} u$ of local functions $u$, while
$L^2_{\text{sol}}(M)$ is the orthogonal complement of
$L^2_{\text{pot}}(M)$ in $L^2 (M)$.
\end{Def}

In Lemmata \ref{piazzola}, \ref{budapest} and \ref{miele} below we
collect some identities relating $\nabla ^{(\o)}$, $\nabla^{(\o),
*}$ and the spatial gradient $\nabla _{\o,e} ^\e$, $e \in \cB$,  defined as follows.
Given a  function $u : \e \cC (\o ) \rightarrow \bbR$, the gradient
$\nabla _{\o,e} ^\e u$ is the function $\nabla _{\o,e} ^\e u : \e
\cC (\o ) \rightarrow \bbR$ defined as
$$
\nabla _{\o,e} ^\e u(x) = \begin{cases} \o (x/\e, x/\e+e ) \frac{
u(x+\e
e)-u(x)}{\e} & \text{ if }  x , x+ \e e  \in \e \cC (\o)\,,\\
0 & \text{ otherwise}\,.
\end{cases}
$$
It can be written as
$$
\nabla_{\o,e}^\e u (x) = \t_{x/\e} \o (0,e) \nabla_e^\e u (x)
$$
where the gradient $\nabla_e^\e u $ is defined as
$$
\nabla _e  ^\e u(x) = \begin{cases}   \frac{ u(x+\e
e)-u(x)}{\e} & \text{ if }  x , x+ \e e  \in \e \cC (\o)\,,\\
0 & \text{ otherwise}\,.
\end{cases}
$$

Lemma \ref{piazzola}  explains why $\nabla ^{(\o)*}$ is called
divergence, or adjoint gradient:
\begin{Le}\label{piazzola}
Given  functions    $u\in L^2 (\mmu) $ and  $v\in L^2(M)$,  it holds
\begin{equation}\label{integrale}
\int_{\O  \times \cB}  v \nabla ^{(\o)} u\, dM= -\int _{\O}  \bigl(
\nabla^{(\o)*} v\bigr)\,u\, d\mu  \,.
\end{equation}
%
In particular,  $\int _{\O  \times \cB} v dM =0$ for any $v\in L^2
_{\text{pot}} (M)$, while    a square integrable form $v\in L^2(M)$
is solenoidal if and only if $\nabla^{(\o)*} v (\o)=0$ for $\mmu$
a.a. $\o$.
\end{Le}
\begin{proof}
By definition
\begin{multline}\label{scherzo}
\int_{\O  \times \cB}  v \nabla ^{(\o)} u\,  dM=\sum _{e\in \cB }
\int _{\O} \mu (d\o )  \o (0,e) v (\o , e ) \bigl( u ( \t_e \o ) - u
(\o) \bigr)\\= - \sum _{e\in \cB} \int _\O \mu ( d \o)  \o (0,e) v
(\o, e) u (\o ) + \sum _{e\in \cB } \int_\O \mu(d\o)\o (0,e) v(\o ,
e) u (\t_e \o )\,.
\end{multline}
Since
$$
\mu (d \o)  \o (0,e) = \bbQ (d\o) \bbI _{0\in \cC (\o)} \o(0,e)
 = \bbQ (d\o) \bbI _{e \in \cC (\o) }  \o (0,e)= \bbQ (d\o)
\bbI  _{0 \in \cC (\t_e \o ) }  \t_e \o  (0, -e),
$$
$v(\o , e)= v\bigl( \t_{-e} (\t_e \o ),e \bigr) $ and $ \bbQ (d\o)=
\bbQ ( d \t_e \o )$, we can conclude that
\begin{multline*}
 \int_\O \mu(d\o) \o
(0,e) v(\o , e) u (\t_e \o )= \int _{\O} \bbQ ( d \t_e\o )\bbI  _{0
\in \cC (\t_e \o ) }  \t_e \o  (0, -e)v\bigl( \t_{-e} (\t_e \o ),e
\bigr)u (\t_e \o )\\= \int _{\O} \bbQ ( d \o )\bbI  _{0 \in \cC (\o
) }  \o (0, -e)v(\t_{-e} \o , e )
 u (\o ) = \int _{\O} \mu ( d \o ) \o (0, -e)v(\t_{-e} \o , e )
 u (\o ) \,.
\end{multline*}
Hence the last sum in (\ref{scherzo}) can be rewritten as
$$
\sum _{e\in \cB } \int_\O \mu(d\o)\o (0,e) v(\o , e) u (\t_e \o ) =
\sum _{e\in \cB } \int _{\O} \mu ( d \o ) \o (0, e)v(\t_{e} \o , -e
)
 u (\o )\,.
 $$
 The above identity and (\ref{scherzo}) allows to conclude the proof
 of (\ref{integrale}), while the second part of the lemma follows easily
 from (\ref{integrale}).
\end{proof}

We point out another integration by parts formula.
\begin{Le}\label{budapest}
Let $e \in \cB$,    $u\in L^1 (\mu^\e _\o)$,   $\psi \in B (\O)$ and
define $v: \O \times \cB \rightarrow \bbR$ as  $$ v
(\o,e'):=\psi(\o) \delta _{e,e'}\,.$$ Then,
\begin{equation}\label{martina}
\int _{\bbR ^d}  \nabla_{\o,e} ^\e u (x) \psi (\t_{x/\e} \o ) \mu
_\o ^\e (dx)=- \e ^{-1} \int _{\bbR^d} u(x) \bigl(\nabla ^{(\o)*} v
 \bigr) ( \t_{x/\e} \o )\mu _\o ^\e (dx)\,.
\end{equation}
\end{Le}
\begin{proof}
By definition of $v$, we have
\begin{equation}\label{ungheria}
\nabla ^{(\o)*}v(\o)=\o(0,e)\psi (\o) -\o(0,-e) \psi (\t_{-e} \o)\,.
\end{equation}
Moreover, since $ \bbI_{z \in \cC (\o) }\o (z,z+e)= \bbI_{z +e \in
\cC (\o) }\o (z,z+e)$, we can write
\begin{multline}\label{zingari}
 \sum _{ z\in \cC (\o)} \o (z,z+e) \left( u(\e z+\e e) -
u(\e z )\right) \psi (\t_z \o) = \\-  \sum _{z\in \cC (\o)} u(\e
z)\left( \o (z,z+e)\psi (\t_z \o) -  \o (z,z-e) \psi ( \t_{z-e}
\o)\right) \,.
\end{multline}
Identities (\ref{ungheria}) and (\ref{zingari})  allow to conclude
the proof of (\ref{martina}).
\end{proof}

Finally, we point out  the simple identities
\begin{align}
& \nabla ^\e_{\o,e} \left( a(x)b(x)\right)=\left( \nabla ^\e_{\o,e}
a(x)\right)b(x+\e e )+ a(x)\left( \nabla^\e_{\o,e}
b(x)\right)\,,\label{sara}\\
& \nabla ^\e_e \left( a(x)b(x)\right)=\left( \nabla ^\e_e
a(x)\right)b(x+\e e )+ a(x)\left( \nabla^\e_e
b(x)\right)\,,\label{sarabis}
\end{align}
valid for all functions $a,b:\e\cC (\o) \rightarrow \bbR^d $. In
what follows, (\ref{sara}) and (\ref{sarabis}) will be frequently
used without explicit mention.

\begin{Le}\label{miele}
Let $u\in L^2 ( \mmu)$.  Suppose that for all  functions $\psi \in
C(\O) $ and for all $e\in \cB$ it holds \begin{equation} \int _\O
u(\o) \nabla ^{(\o)*} v(\o ) \mmu (d\o)=0 \,,\qquad v
(\o,e'):=\psi(\o) \delta _{e,e'}\,.
\end{equation} Then,
$u$ is constant $\mu$--almost everywhere.
\end{Le}
\begin{proof}
Due to Lemma \ref{piazzola}, for all  functions $\psi\in C(\O) $ and
for all $e\in \cB$ it holds
\begin{equation}
0= \int_{\O\times \cB }\bigl( \nabla ^{(\o)} u \bigr) v \, dM= \int
_\O \mu (d \o) \o (0,e) \left( u (\t_e\o) - u(\o) \right) \psi
(\o)\,.
\end{equation}
Hence, \begin{equation} \label{nerone}
\bbI _{\{0,e\} \in \cE (\o)  }\left( u (\t_e\o) - u(\o) \right)
=0\,, \qquad \bbQ\text{--}a.s.
\end{equation}
Due the translation invariance of $\bbQ$ we conclude that
\begin{equation} \label{neronebis}
\bbI _{\{x,x+e\} \in \cE (\o)  }\left( u (\t_{x+e}\o) - u(\t_x \o)
\right) =0\,, \qquad \forall x \in \bbZ^d\,, \;\bbQ\text{--}a.s.
\end{equation}
 Since $\cC(\o)$ is connected, (\ref{neronebis}) is
equivalent to say that  for $\bbQ$--a.a. $\o$ there exists a
constant $a(\o) $ such that $u(\t_x\o)=a(\o)$ for all $x \in
\cC(\o)$. Trivially, the function $a(\o)$ is translation invariant.
Hence, due to the ergodicity of $\bbQ$ we can conclude that $a(\o)$
is constant $\bbQ$--a.s. Since $a(\o)= u(\o)$ if $0 \in \cC(\o)$, we
conclude that $u(\o)$ is constant for $\mu$--a.a. $\o$.
\end{proof}

\medskip

We have now all the tools in order to define the set of regular
environments $\O_*$. To this aim we first observe that $L^2
_\text{sol}$  is separable since it is  a subset of the separable
metric space $L^2(M)$. We fix once and for all a sequence $\{\psi
_j\}_{j\geq 1}$  dense in $L^2_{\text{sol}}$. Since elements of
$L^2(M)$ are equivalent if, as functions, they differ on a zero
measure set, we fix a representative $\psi_j$ and from now on we
think of $\psi_j$  as pointwise function $\psi_j:\O\times \cB
\rightarrow \bbR$. Since $\psi_j \in L^2 (M)$ it must be
\begin{equation}\label{salerno}
\Psi_{j,e}(\o):=\sqrt{\o(0,e) } \psi _j(\cdot, e) \in L^2 (\mu)\,.
\end{equation}
For each $j\geq 1$ and $e \in \cB$ we fix a sequence of continuous
functions $f_{j,e}^{(k)} \in C(\O)$ such that $f_{j,e}^{(k)}$
converges to $ \Psi_{j,e} $ in $L^2 (\mu)$ as $k\rightarrow \infty$.
We first make a simple observation:
\begin{Le}\label{zorro}
Given  $j\geq 1$, define  the Borel set $\O_{*,j}$  as the set of
configurations $\o \in \O_0$ such that
\begin{align}
& \nabla ^{(\o)*}\psi_j (\t_x\o)=0  \;\; \forall x \in \cC (\o)\,,\label{trieste1}\\
& \lim _{\e\downarrow 0} \int_{[-n,n]^d} \bigl( f^{(k)}_{j,e}
(\t_{z/\e}\o ) -\Psi_{j,e} ( \t_{z/\e} \o ) \bigr)^2\mu^\e_\o(dz)
=(2n)^d  \bigl\|f^{(k)}_{j,e}-\Psi_{j,e}
\bigr\|^2_{L^2(\mu)}\label{trieste2}
\end{align}
for each $ e \in \cB$ and $k,n \geq 1$. Then $\bbQ( \O_{*,j})=1$.
\end{Le}
\begin{proof}
Let us define the set $A_x$ as
$$ A_x :=\{ \o \in \O_0\,:\,  \nabla ^{(\o)*}\psi_j (\t_x\o)\not =0 \text{ and  } x \in \cC (\o) \}\,,
\qquad x \in \bbZ^d\,.$$ Due to Lemma \ref{piazzola}, $\bbQ(A_0)=\mu
(A_0)=0 $. Since $\o \in A_x $ if and only if $\t_x \o \in A_0$, by
the translation invariance of $\bbQ$ we obtain that $\bbQ(A_x)=0$
for all $x \in \bbZ^d$. Hence, setting $A= \cup _{x\in \bbZ^d} A_x$,
it must be $\bbQ(A)=0$. Since $\O_0\setminus A$ coincides with the
set of $\o \in \O_0$ satisfying (\ref{trieste1}), we only need to
prove that (\ref{trieste2}) is satisfied $\bbQ$--a.s. for each
$k,n\geq 1 $ and $e \in \cB$. This is a direct consequence of the
$L^1$--ergodic theorem.
\end{proof}
\medskip
Recall the definition of $\O_1 \subset \O_0$ given before Theorem
\ref{prop_hom}. We can finally define the set $\O_*$:
\begin{Def} We define the  set $\O_*$ of regular environments as
$$ \O_*:= \O_1\cap \left(\cap _{j=1}^\infty \O_{*,j}\right) \subset \O_0\,.$$
\end{Def}

\bigskip

 We conclude this section by reformulating   the
variational characterization (\ref{varcar}) of the diffusion matrix
$\cD$ in terms of square integrable forms. To this aim,  given a
vector $\xi \in \bbR^{\cB_*} $, we write $w^\xi$ for the square
integrable form
\begin{equation}\label{elefante}
w^\xi (\o, \pm e) := \pm \xi _e \,,\qquad  e\in \cB_*, \;\;
\o\in\O\,.
\end{equation}
Let  $\p: L^2 (M) \rightarrow L^2 _{\text{sol}}(M) $  be the
orthogonal projection of $L^2(M)$ onto  $L^2 _{\text{sol}}(M) $ and
let  $\Phi$ be  the bilinear form on $\bbR ^{\cB_*}\times
\bbR^{\cB_*}$ defined as
$$
\Phi (\zeta, \xi) = ( w^\zeta, \p w^\xi )_{L^2(M)}\,,
$$
where $( \cdot , \cdot  )_{L^2(M)}$ denotes the inner product in
$L^2 (M)$.

 Since $\Phi$ is bilinear and symmetric, there exists a
symmetric matrix $D$ indexed on $\cB_*\times \cB_*$ such that
\begin{equation}\label{zenfira}
( \zeta, D\xi )=  ( w^\zeta, \p w^\xi )_{L^2(M)}\,.
\end{equation}
We  give an integral representation of $D\xi$ which will be useful
in what follows. Since
\begin{equation}\label{totti}
 ( w^\zeta, \p w^\xi )_{L^2(M)}
=\sum _{e\in \cB_*} \zeta_e \int _\O \mmu (d\o)\left[
 \o (0,e) \bigl(\p w ^\xi  \bigr) (\o, e) -\o (0,-e) \bigl(\p w ^\xi  \bigr) (\o, -e)\right]\,,
 \end{equation}
 it must be
 \begin{equation}\label{anni}
( D \xi)_e = \int _\O \mmu (d\o)\left[
 \o (0,e) \bigl(\p w ^\xi  \bigr) (\o, e) -\o (0,-e) \bigl(\p w ^\xi  \bigr) (\o, -e)\right]\,.
  \end{equation}
Moreover, due to the definition of orthogonal projection, we get
\begin{equation}\label{maria0}
\begin{split}
(\xi , D\xi )& =  ( w^\xi, \p w^\xi )_{L^2(M)} = \| \p w^\xi \| _{L^2(M)}^2 =
\inf _{v \in  L^2 _{\text{pot}}  (M)}\|  w^\xi-v  \| _{L^2(M)}^2\\
&  = \inf _{\psi \in B(\O)   }\|  w^\xi-\nabla ^{(\o) }\psi \|
_{L^2(M)}^2\,.
\end{split}
\end{equation}
By definition,
\begin{equation}\label{maria}
\begin{split}
 \|  w^\xi-\nabla ^{(\o) }\psi  \| _{L^2(M)}^2  &=   \sum _{e\in
\cB_*} \int _\O  \mu  (d\o) \o (0,e )  \bigl[\xi _e -\hat
\o(0,e)\bigl( \psi
(\t_e \o)-\psi(\o) \bigr)\bigr] ^2\\
&  +
 \sum _{e\in \cB_*} \int _\O  \mu  (d\o) \o(0,-e) \bigl[ -\xi _e -\hat \o (0, -e)  \bigl(  \psi (\t_{-e} \o)-\psi(\o)
\bigr)\bigr] ^2  \\
&   = \sum _{e\in \cB_*}  \int _\O \mu (d\o) \o (0,e)  \bigl( \xi _e
- \psi (\t_e \o)+\psi(\o)
\bigr) ^2    \\
& +  \sum _{e\in \cB_*} \int _\O  \mu  (d\o) \o (0,-e)
  \bigl( -\xi _e - \psi (\t_{-e} \o)+\psi(\o)
\bigr) ^2\,.
 \end{split}
\end{equation} We can rewrite the last term in a more useful form.
In fact,  due to the translation invariance of $\bbQ$, we get
\begin{equation}\label{filippa}
\begin{split}
&\int _\O  \mu  (d\o) \o (0,-e)
  \bigl( -\xi _e - \psi (\t_{-e} \o)+\psi(\o)
\bigr) ^2=\\
& \int _\O  \bbQ (d\o) \bbI_{0,e\in \cC (\t_{-e} \o) } \t_{-e} \o
(0,e)  \bigl( -\xi _e - \psi (\t_{-e} \o) +\psi\bigl(\t_e
(\t_{-e}\o))
\bigr) ^2=\\
&  \int _\O  \bbQ (d\o) \bbI_{0,e\in \cC (\o) }\o (0,e)  \bigl( \xi
_e +\psi(\o)   - \psi (\t_e \o)\bigr) ^2 =
 \int _\O  \mu  (d\o)\o (0,e)  \bigl( \xi
_e - \psi (\t_e \o)+\psi(\o) \bigr) ^2 \,.
\end{split}
\end{equation}
 Due to (\ref{maria0}), (\ref{maria}) and (\ref{filippa}) we
 conclude that
 \begin{equation}\label{mozart111}
 (\xi , D\xi )= \inf_{\psi\in B (\O) }2
 \sum _{e\in \cB_*} \int _\O  \mu (d\o ) \o (0,e)  \bigl( \xi _e - \psi (\t_e \o)+\psi(\o)
\bigr) ^2 \,.
\end{equation}
 In particular, the matrix $D$ is related to the matrix $\cD$ via
 the identity
 \begin{equation}\label{sensi}
 D=2m \mathcal{D}\,.
 \end{equation}

From the above observations and the  non--degeneracy of $\cD$ given
by hypothesis (H3) we get:
  \begin{Le}\label{ranatan1000}
  The vectorial space given by the vectors
  $$
 \left( \int _\O \mmu (d\o)\left[
 \o (0,e) \psi (\o, e) -\o (0,-e) \psi (\o, -e)\right]
 \right)_{e\in \cB_*} ,\qquad \psi \in L^2 _{\text{sol}}\,,
 $$
 coincides with $\bbR^{\cB_*}$.
   \end{Le}
   \begin{proof}
   If the statement was not true, then there would exist $\xi \in \bbR^{\cB_*}\setminus\{0\}$ such that
   $$
   \sum _{e\in \cB_*}\xi _e   \int _\O \mmu (d\o)\left[
 \o (0,e) \psi (\o, e) -\o (0,-e) \psi (\o, -e)\right]
 =0, \qquad \forall \psi \in L^2 _{\text{sol}}\,. $$
 In particular, the above identity would hold  with $\psi = \p w ^\xi $. Due to (\ref{zenfira}) and  (\ref{totti}), this would imply that
 $
 (\xi, D \xi )=0$, which is absurd due to hypothesis (H3).
   \end{proof}
Finally, we conclude  with   a simple but crucial observation.
   Given $\xi \in \bbR^{\cB_*}$, there exists a unique form $v \in L^2 _{\text{pot}}$ such that
     $w^\xi +v\in L^2_{\text{sol}}$. In fact, these requirements imply that
      $w^\xi +v=\p w^\xi$.

\section{Two--scale convergence}\label{lampione}

In this section we  analyze the  weak two--scale convergence for our
disordered model. We recall that $\O_*$ denotes the set of regular
environments $\o$ defined in the previous section, and we recall
that $(\cdot, \cdot)_{\mu ^\e _\o}$ and $\|\cdot \|_{\mu ^\e_\o} $
denote respectively   the inner product and the norm in $L ^2 (\mu
^\e_\o)$. In our context the two--scale convergence
\cite{ZP}[Section 5] can be defined as follows:
\begin{Def}\label{landim}Fix $\o \in \O_*$.
Let $v ^\e$  be  a family of functions parameterized by $\e>0$ such
that  $v ^\e \in L^2 (\mu ^\e _\o)$. Then the function $v \in L^2
(\bbR^d \times \O, dx \times \mmu )$
 is  the weak two--scale limit of $ v ^\e $ as
$\e\downarrow 0$ (shortly, $v^\e\stackrel{2}{\rightharpoonup}v$) if
the following two conditions are fulfilled:
\begin{equation}\label{limitato}
\limsup_{\e\downarrow 0} \| v ^\e \|_{\mu ^\e _\o }
<\infty\,,
\end{equation}
and
\begin{equation}\label{limitone}
\lim _{\e\downarrow 0} \int _ {\bbR^d} v^\e  (x) \varphi (x) \psi
\bigl( \t_{x/\e} \o\bigr) \mu _\o^\e (dx)=\int _{\bbR^d}dx  \int
_{\O} v(x, \o' ) \varphi (x) \psi (\o ') \mmu(d\o')  \,,
\end{equation}
for all $\varphi \in C^\infty_c (\bbR^d)$ and  $\psi \in C(\O)$.
\end{Def}

Let us first collect some technical results concerning the weak
two--scale convergence. For the next lemma, recall the definition of
the function $\Psi_{j,e}\in L^2(\mu)$ given in \eqref{salerno}.
\begin{Le}\label{repubblica}
Fix $\o \in \O_*$ and suppose that $L^2(\mu^\e_\o)\ni
v^\e\stackrel{2}{\rightharpoonup}v\in L^2(\bbR^d \times \O,dx\times
\mu)$.
 Then, for each $j\geq 1$, $e \in \cB$, $\psi \in C(\O)$ and
 $\varphi \in C_c (\bbR^d)$,
it holds
\begin{equation}\label{limitonebis}
\lim _{\e\downarrow 0} \int _ {\bbR^d} v^\e  (x) \varphi (x)  \psi
\bigl( \t_{x/\e} \o \bigr) \Psi_{j,e} \bigl( \t_{x/\e} \o \bigr) \mu
_\o^\e (dx)=\int _{\bbR^d}dx \int _{\O} v(x, \o' ) \varphi (x) \psi
(\o ') \Psi_{j,e}(\o') \mmu(d\o')  \,.
\end{equation}
\end{Le}
\begin{proof}
Suppose that the support of $\varphi$ is included in $[-n,n]^d$.
Recall the definition of the functions $f^{(k)}_{j,e}\in C(\O)$
given in Section \ref{picche}.  Then, by Schwarz inequality, we get
\begin{multline}
\left|\int _ {\bbR^d} v^\e  (x) \varphi (x) \psi( \t_{x/\e} \o)
\bigl[\Psi_{j,e} \bigl( \t_{x/\e} \o \bigr)- f^{(k)}_{j,e} \bigl(
\t_{x/\e} \o \bigr)\bigr]
\mu _\o^\e (dx)\right|\leq \\
 \|\varphi\|_\infty \|\psi\|_\infty \| v^\e
\|_{\mu^\e_\o} \Big\{ \int_{[-n,n]^d} \bigl( \Psi_{j,e} (\t_{x/\e}\o
) -f^{(k)}_{j,e} ( \t_{x/\e} \o)
\bigr)^2\mu^\e_\o(dx)\Big\}^{1/2}\,.
\end{multline}
Since $\o \in \O_*\subset \O_{*,j}$ and since
$f^{(k)}_{j,e}\rightarrow \Psi_{j,e}$ in $L^2(\mu)$, we conclude
that the upper limit of the r.h.s.   as $\e\downarrow 0$ and then
$k\uparrow \infty$ is zero. On the other hand, since
$v^\e\stackrel{2}{\rightharpoonup}v$ and $f^{(k)}_{j,e}\rightarrow
\Psi_{j,e}$ in $L^2(\mu)$, we obtain that
\begin{multline}
\lim_{k\uparrow \infty }\lim_{\e\downarrow 0} \int _ {\bbR^d} v^\e
(x) \varphi (x) \psi (\t_{x/\e} \o) f^{(k)}_{j,e} \bigl( \t_{x/\e}
\o \bigr) \mu _\o^\e (dx)= \\\lim_{k\uparrow \infty}\int _{\bbR^d}dx
\int _{\O} v(x, \o' ) \varphi (x)\psi(\o') f^{(k)}_{j,e} (\o ')
\mmu(d\o')=\\
 \int _{\bbR^d}dx \int
_{\O} v(x, \o' ) \varphi (x) \psi(\o')\Psi_{j,e} (\o ')
\mmu(d\o')\,.
\end{multline}
This allows to get (\ref{limitonebis}).
\end{proof}

By the same arguments leading to \cite{ZP}[Lemma 5.1] and
\cite{Z}[Prop. 2.2] one can easily prove the following result:
\begin{Le}\label{lemmino}
Fix $\o \in \O_*$. Suppose that the  family of functions $v^\e \in
L^2 (\mu^\e_\o)$ satisfies (\ref{limitato}). Then from each sequence
$\e_k $ converging to zero, one can extract a subsequence $\e_{k_n}$
such that $v^\e$ converges along $\e_{k_n}$  to some $v \in L^2
(\bbR^d \times \O, dx \times \mmu )$ in the sense of weak two--scale
convergence.
\end{Le}
We give the proof for the reader's convenience:

\begin{proof}
Given $ \varphi \in C_c ^\infty (\bbR ^d)$ and   $\psi \in C(\O)$,
we can bound
\begin{multline*}
 \limsup _{k \rightarrow \infty} \left| \int_{\bbR ^d} v^{\e_k}
(x)  \varphi (x) \psi (\t_{x/\e_k} \o ) \mu^{\e_k} _\o (dx)
\right|\leq \\
 \limsup _{k \rightarrow \infty} \| v^{\e_k}\| _{\mu ^{\e_k} _\o}
\left( \int_{\bbR ^d}  \varphi^2  (x) \psi^2 (\t_{x/\e_k} \o )
\mu^{\e_k} _\o (dx) \right)^{1/2}\leq \\
C(\o)  \left(\int_{\bbR^d} \varphi ^2(x) \, dx \int_{\O} \psi^2 (\o'
) \mu (d\o') \right)^{1/2}= C(\o)  \bigl \| \varphi  \psi \bigr
\|_{L^2( \bbR ^d \times \O , dx \times \mu)}
\end{multline*}
(note that the first estimate follows from Schwarz inequality, while
the second one follows from (\ref{presilla}) and (\ref{limitato})).

Using a standard diagonal argument and the separability of the space
of test functions  $\varphi,\psi$, we can conclude that there exists
a subsequence $\{\e_{k_n}\} _{n\geq 1}$ along which  the limit in
the l.h.s. of (\ref{limitone}) exists and can be extended to a
continuous linear functional on $L^2(\bbR ^d \times \O, dx \times
\mu)$. Therefore, this limit can be written as the inner product in
$L^2(\bbR ^d \times \O, dx \times \mu)$ with  a suitable function
$v$.

\end{proof}

\smallskip

In what follows, we will apply the concept of weak two--scale
convergence to the solution $u^\e_\o\in L^2 (\mu^\e_\o)$ of
(\ref{mammina}) and to its gradients, for a fixed sequence $f^\e_\o
\strano f$.
 To this aim we start with some
simple observations.

 We note that, given
  $u,v\in L^2 (\mu^\e _\o)
$, it holds
\begin{multline}\label{rinog}
(u, -\cL ^\e _\o v )_{\mu ^\e _\o} = \frac{\e^{d-2} }{2} \sum_{z\in
\cC (\o) } \sum _{e \in \cB } \o (z,z+e) \left[ u(\e z + \e e) - u
(\e z ) \right]\left[ v(\e z + \e e) - v (\e z ) \right]
 \\= \frac{\e^d }{2} \sum_{z\in \cC (\o) } \sum _{e \in \cB }
\o (z,z+e) \nabla _e^\e  u (\e z ) \nabla ^\e_e  v( \e z ) =
\frac{1}{2} \sum _{e \in \cB } \int _{\bbR^d} \mu ^\e _\o (dx)
\t_{x/\e} \o   (0,e) \nabla ^\e_e u (x) \nabla _e^\e v( x) \,.
\end{multline}
In particular, we can write
\begin{multline}\label{diri}
(u^\e_\o,-\cL ^\e_\o u^\e_\o)_{\mu^\e_\o} = \frac{1}{2}\sum_{e\in
\cB}
 \int _{\bbR^d} \mu ^\e _\o (dx) \t_{x/\e} \o   (0,e)
\bigl( \nabla ^\e_e  u ^\e _\o (x) \bigr)^2 = \\
\frac{1}{2}\sum_{e\in \cB} \Big\| \sqrt{\t_{x/\e} \o (0,e) }
\nabla^\e_e u^\e_\o (x) \Big\|_{\mu^\e_\o}^2 \,.
\end{multline}
Moreover, taking the inner product of (\ref{mammina}) with
$u^\e_\o$, we obtain
$$
\l \| u^\e_\o\|_{\mu^\e_\o}^2\leq \l (u^\e_\o, u^\e_\o)_{\mu^\e_\o}
+ (u^\e_\o,-\cL ^\e_\o u^\e_\o)_{\mu^\e_\o}= (u^\e_\o, f^\e _\o
)_{\mu^\e_\o}\leq
\| u^\e_\o\|_{\mu^\e_\o} \| f^\e_\o  \|_{\mu^\e_\o}\,.
$$
Hence, since   $f^\e_\o\strano f$,
 for any $\l>0$ it holds that
\begin{equation}\label{mazinga}
\sup _{\e>0} \| u^\e_\o\|_{\mu^\e_\o}^2<\infty\,,\;\;\sup _{\e>0} \,
(u^\e_\o,-\cL ^\e_\o u^\e_\o)_{\mu^\e_\o}<\infty\,, \;\; \sup_{\e>0,
e\in \cB } \Big\| \sqrt{\t_{x/\e} \o (0,e) } \nabla^\e_e u^\e_\o (x)
\Big\|_{\mu^\e_\o}<\infty
 \,.
\end{equation}

\begin{Le}\label{lago}
Fix $\tilde \o \in \O_*$. The family $u^\e_{\tilde \o}$ converges
along a subsequence to a function $u^0\in L^2 (\bbR^d \times \O, dx
\times \mmu )$ in the sense of weak two--scale convergence and $u^0$
does not depend on $\o$, i.e. $u^0\in L^2 (\bbR^d, dx)$.
\end{Le}
\begin{proof}
Due to Lemma \ref{lemmino}, the sequence $u^\e_{\tilde \o}$
converges along a subsequence $\e_k\downarrow 0 $ to a function
$u^0\in L^2 (\bbR^d \times \O, dx \times \mmu )$ in the sense of
weak two--scale convergence. In order to simplify the notation, we
suppose that this convergence holds for $\e\downarrow 0$. We need to
prove that $u^0$ does not depend on $\o$. To this aim, fix $e\in
\cB$, $\varphi \in C_c ^\infty (\bbR ^d)$ and a  function $\psi \in
C(\O)$. We define $v(\o,e')=\psi (\o) \d_{e,e'}$. Due to the
definition of weak two--scale convergence, it holds
\begin{equation}\label{gianni}
\lim _{\e \downarrow 0 } \int _{\bbR ^d} u _{\tilde \o}^\e (x)
\varphi(x) \nabla ^{(\o)*} v ( \t _{x/\e} {\tilde \o} )\mu ^\e
_{\tilde\o} (dx) = \int_{\bbR^d} dx \int _\O \mmu (d\o) u^0(x,\o)
\varphi (x) \nabla ^{(\o)*} v (\o )\,,
\end{equation}
while due to Lemma \ref{budapest}
\begin{equation}\label{talamanca}
 \int _{\bbR ^d} u _\oo^\e (x) \varphi(x)
\nabla ^{(\o)*} v ( \t _{x/\e} \oo )\mu ^\e _\oo (dx)= - \e \int
_{\bbR ^d}  \nabla_{\o,e} ^\e \left( u^\e_\oo (x) \varphi (x)
\right) \psi ( \t _{x/\e} \oo )\mu ^\e _\oo (dx)\,.
\end{equation}
The r.h.s. in (\ref{talamanca}) is bounded by
\begin{equation}
\e\|\psi\|_\infty \int _{\bbR ^d}\left|\nabla_{\o,e} ^\e \left(
u^\e_\oo(x) \varphi (x) \right)\right| \mu ^\e _\oo (dx) \leq I_1
+I_2\,,
\end{equation}
where
\begin{align*}
& I_1=\e\|\psi\|_\infty \int _{\bbR ^d} \left|\nabla_{\o,e} ^\e
u^\e_\oo (x)\cdot  \varphi (x+\e e )\right|\mu ^\e _\oo (dx)\, ,\\
& I_2=\e\|\psi\|_\infty \int _{\bbR ^d} \left| u^\e_\oo (x)\cdot
\nabla_{\o,e} ^\e \varphi (x )\right|\mu ^\e _\oo (dx)\, .
\end{align*}
By Schwarz inequality and (\ref{diri})  we can bound
\begin{equation}
\begin{split}
I_1 & \leq \e \| \psi \|_\infty \left[ \int_{\bbR^d} \t_{x/\e} \oo
(0,e) \, \bigl ( \nabla _e ^\e u^\e_\oo (x) \bigr) ^2 \mu^\e_\oo
(dx) \right]^{1/2} \left[ \int_{\bbR^d}\t_{x/\e} \oo(0,e)\, \varphi
(x+\e e )^2 \mu^\e_\oo (dx) \right]^{1/2} \\& \leq \e\, c(\psi,
\varphi) \bigl( u^\e_\oo , - \cL ^\e _\oo u ^\e _\oo\bigr)^{1/2}
_{\mu^\e _\oo }\,,
\end{split}
\end{equation}
for a suitable positive constant $c(\psi,\varphi)$ depending on
$\psi$ and $\varphi$. Due to (\ref{mazinga}),  we obtain that $I_1
\leq c(\psi,\varphi,\oo) \e $.

Moreover, by Schwarz inequality we have
$$ I_2  \leq \e \|\psi\|_\infty \|  u^\e _\oo  \|_{\mu_\oo^\e} \left
\|\nabla ^\e _{\o,e} \varphi \right\|_{\mu^\e_\oo}\,,
$$
and again from (\ref{mazinga}) we deduce that $I_2\leq
c(\psi,\varphi,\oo) \e $.
  Hence the r.h.s. of
(\ref{talamanca}) is bounded by $c\e$ and due to (\ref{gianni}) we
get that $$ \int_{\bbR^d} dx \int _\O \mmu (d\o) u^0(x,\o) \varphi
(x) \nabla ^{(\o)*} v (\o )=0\,.$$ Since this holds for all $\varphi
\in C^\infty_c (\bbR^d)$ we get   that
$$
\int _\O \mmu (d\o) u^0(x,\o)  \nabla ^{(\o)*} v (\o )=0
$$
for Lebesgue a.a. $x\in \bbR^d$. Due to separability, we conclude
that  for Lebesgue a.a. $x\in \bbR^d$ the above identity is valid
for all  $v$ of the form $v(\o, e')=\psi (\o) \delta_{e,e'}$, for
some  function $\psi\in C(\O)$ and some $e\in \cB$.  By Lemma
\ref{miele} we conclude that for these points $x$, the function
$u^0(x, \cdot)$ is constant $\mmu$--almost everywhere. This
concludes the proof.
\end{proof}

In what follows, $u^0$ will be as in Lemma \ref{lago} for a fixed
$\oo\in \O_*$. We will prove at the end that $u^0$ coincides with
the solution of (\ref{caldo}) and in particular that $u^0$ does not
depend on $\oo$.

\begin{Le}\label{bracciano} Fix $\oo\in \O_*$.
The function $u^0$ belongs to the Sobolev space $H^1(\bbR^d,  dx )$.
Moreover, along  a suitable subsequence and for all $e\in \cB$ it
holds
\begin{align}
& u ^\e_\oo (x)\stackrel{2}{ \rightharpoonup} u^0(x)\,,\label{mora}\\
& \sqrt{\t_{x/\e}   \oo (0,e) } \nabla ^\e_e u^\e_\oo(x)
\stackrel{2}{\rightharpoonup} v^0_e(x, \o) \label{mela}
\end{align}
for some $v^0_e \in L^2 ( \bbR ^d \times \O , dx \times \mu ) $.

Given $x\in \bbR^d$, consider the forms $\theta_x$, $\Gamma_x$ in
$L^2(M)$  defined as
\begin{align}
& \theta_x : \O \times \cB \ni (\o,e) \rightarrow v^0_e
(x,\o)/\sqrt{\o(0,e) } \in \bbR\,,\label{mamma1}\\
& \G_x : \O \times \cB \ni (\o,e) \rightarrow \partial_e u^0 (x)\in
\bbR^d \,,\label{mamma2}
\end{align}
where $ \partial_e u^0 (x)$ denotes a representative of the weak
derivative in $L^2 (dx)$ of $u^0$,  along the  direction $e$. Then,
for Lebesgue a.a. $x \in \bbR^d$, it holds
\begin{equation}\label{seppioline}
\theta_x \in L^2_{\text{sol} } (M)\,, \qquad \theta _x = \p \G_x \,,
\end{equation}
where $\p : L^2(M)\rightarrow L^2_{\text{sol}} (M)$ is the
orthogonal projection onto $L^2_{\text{sol}} (M)$.

  \end{Le}
  Note that the form $\theta_x$ is well defined, since for $M$--a.a.
  $(\o, e)$ it holds $\o (0,e)>0$. Moreover, $\theta_x \in L^2(M)$
  for Lebesgue a.a. $x \in\bbR^d$. In fact,
  $$
  \| \theta _x \|_{L^2(M)} ^2 = \sum _{e\in \cB} \int_{\O} \mu (d \o
  ) \o (0,e) \theta_x (\o ,e)^2= \sum  _{e\in \cB} \int_{\O} \mu (d \o
  )  v^0_e (x,\o) ^2  $$
and $v^0_e \in L^2 ( \bbR ^d \times \O , dx \times \mu ) $.

\begin{proof}
(\ref{mora}) follows from Lemma \ref{lago}. At cost to take a
sub--subsequence, due to Lemma \ref{lemmino} and (\ref{mazinga}),
(\ref{mela}) holds for all $e\in \cB$.
\smallskip

Let us prove that $u^0 \in H^1 (\bbR^d, dx)$. To this aim for each
$j\geq 1$ we consider the function $\psi_j:\O\times \cB \rightarrow
\bbR$ introduced before Lemma \ref{zorro} (we recall that $\psi_j
\in L^2_{\text{sol}}(M)$) and we take  a function $\varphi\in
C^\infty_c (\bbR^d)$.
 Then by (\ref{sara}) we can write
\begin{equation}\label{butta}
 \begin{split} \sum _{e\in \cB} \int _{\bbR^d} \bigl(\nabla_{\o,e} ^\e u^\e_\oo(x)
 \bigr)\, &
  \varphi(x)\, \psi_j (\t_{x/\e}\oo,e) \mu^\e_\oo (dx)=\\
&  \sum _{e\in \cB} \int _{\bbR^d} \bigl(\nabla_{\o,e} ^\e
u^\e_\oo(x) \bigr)\,
  \varphi(x+\e e )\, \psi_j (\t_{x/\e}\oo,e) \mu^\e_\oo (dx)+ o(1) =\\ &
 \sum _{e\in \cB} \int _{\bbR^d} \nabla_{\o,e} ^\e \bigl( u^\e_\oo(x)\, \varphi(x)\bigr)
\, \psi_j (\t_{x/\e}\oo,e) \mu^\e_\oo (dx)-
 \\ &
\sum _{e\in \cB} \int _{\bbR^d} u^\e_\oo(x) \bigl(\nabla_{\o,e} ^\e
\varphi(x)\bigr) \, \psi_j (\t_{x/\e}\oo,e) \mu^\e_\oo (dx)+o(1)\,.
\end{split}
\end{equation}
Due to Lemma \ref{budapest}, we can rewrite the first addendum in
the r.h.s. as
\begin{multline}\label{patti}
\sum _{e\in \cB} \int _{\bbR^d} \nabla_{\o,e} ^\e \bigl(
u^\e_\oo(x)\, \varphi(x)\bigr) \, \psi_j (\t_{x/\e}\oo,e) \mu^\e_\oo
(dx)= \\
-\e^{-1} \sum_{e\in \cB} \int _{\bbR^d} u^\e_\oo (x) \varphi (x)
\bigl( \nabla^{(\o)*}\psi_j \bigr)(\t_{x/\e}\oo,e) \mu^\e_\oo
(dx)\,.
\end{multline}
Since $\tilde \o \in \O_*\subset \O_{*,j}$, the r.h.s. is zero (see
(\ref{trieste1})). We conclude that
\begin{multline}
 \sum _{e\in \cB} \int _{\bbR^d} \bigl(\nabla_{\o,e} ^\e u^\e_\oo(x) \bigr)\,
  \varphi(x)\, \psi_j (\t_{x/\e}\oo,e) \mu^\e_\oo (dx)=\\
  -\sum _{e\in \cB} \int _{\bbR^d} u^\e_\oo(x) \partial_e \varphi(x)
\, \psi_j (\t_{x/\e}\oo,e)\t_{x/\e}\oo(0,e) \mu^\e_\oo (dx)+o(1)\,.
\end{multline}
We now take the limit  $\e \downarrow 0$  along the suitable
subsequence of the above indentity.  Since
\begin{equation}\label{daidai}
\nabla_{\o,e} ^\e u^\e_\oo(x) = \sqrt{\t_{x/\e}   \oo (0,e) }\left[
\sqrt{\t_{x/\e}   \oo (0,e) } \nabla ^\e_e u^\e_\oo(x)\right]\,,
\end{equation}  the above identity and Lemma \ref{repubblica} imply
 that
\begin{equation}\label{lunatico}
\begin{split}
& \sum _{e\in \cB}  \int _{\bbR^d} dx \int _\O \mmu (d\o)\sqrt{\o(0,e)} v^0_e (x,\o) \varphi(x) \psi_j (\o,e) =\\
& -
\sum _{e\in \cB}  \int _{\bbR^d} dx\, u^0 ( x) \partial_e \varphi(x) \int _\O  \mmu(d\o)\psi_j (\o, e)\o(0,e) =\\
& -\sum _{e\in \cB_*}  \int _{\bbR^d} dx\, u^0 ( x) \partial_e
\varphi(x)
 \int _\O  \mmu(d\o)\left[\psi_j (\o, e)\o(0,e) - \psi_j (\o,
 -e)\o(0,-e)\right]\,.
 \end{split}
 \end{equation}
Given $\psi \in L^2_{\text{sol}}$, we  define $a(\psi)\in \bbR ^d$
as
\begin{equation*}
a(\psi)_e =  \int _\O  \mmu(d\o)\left[\psi (\o, e)\o(0,e) - \psi
(\o, -e)\o(0,-e)\right]\,, \qquad e\in \cB_*\,.
\end{equation*}
Due to Lemma \ref{ranatan1000},
$\left\{ a(\psi )\,:\, \psi \in L^2 _{\text{sol}} (M)
\right\}=\bbR^{\cB_*}$.
On the other hand, since the map $L^2_\text{sol}\ni \psi \rightarrow
a(\psi) \in \bbR^{\cB_*}$ is continuous, we conclude that
$\{a(\psi_j)\}_{j\geq 1}$ is dense in $\bbR^{\cB_*} $.
 Due to (\ref{lunatico}),
\begin{equation}\label{neutro}
 -\sum _{e\in \cB_*}  a(\psi_j)_e  \int _{\bbR^d} dx\, u^0 ( x) \partial_e \varphi(x) =
 \sum  _{e\in \cB_*} \int _{\bbR^d} dx\, \varphi (x)  f_e (\psi_j, x)
 \end{equation}
 where
 \begin{equation}\label{acido}
 f_e (\psi_j,\cdot )= \int _\O \mmu (d\o)\left(\sqrt{\o(0,e)}  v^0_e (\cdot ,\o)  \psi _j(\o,e)+
 \sqrt{\o(0,-e)}  v^0_{-e} (\cdot ,\o)  \psi _j(\o,-e)\right) \in L^2 (\bbR ^d )\,.
 \end{equation}
As consequence of the density of  $\{a(\psi_j)\}_{j\geq 1}$  in
$\bbR^{\cB_*} $, (\ref{neutro}) and (\ref{acido}), it must be $u^0
\in H^1(\bbR^d, dx )$.

\smallskip

Let us now  prove   (\ref{seppioline}). Since $u^0 \in H^1(\bbR^d,
dx )$,
  we are allowed to rewrite
the first identity in (\ref{lunatico}) as
\begin{equation}\label{lunaticobis}
  \sum _{e\in \cB}    \int _{\bbR^d} dx  \, \varphi(x)
 \int _\O \mmu (d\o)\left( \sqrt{\o(0,e)} v^0_e (x,\o)-   \o (0,e)\partial _e u^0 (x)      \right)   \psi_j (\o,e)=0\,.
\end{equation}
Then,
 by means of the arbitrariness  of
$\varphi$ and separability arguments, we get that for Lebesgue a.a.
$x\in \bbR^d$
\begin{multline*}
  \sum _{e\in
\cB} \int _\O \mmu (d\o)\left(\sqrt{\o(0,e)}  v^0_e (x,\o)- \o(0,e)
\partial _e
u^0 (x) \right) \psi_j (\o,e)=\\
 \sum _{e\in
\cB} \int _\O \mmu (d\o) \o (0,e) \left( \frac{ v^0_e
(x,\o)}{\sqrt{\o(0,e)}}   -
\partial _e u^0 (x) \right) \psi_j (\o,e)=  0\,, \qquad \forall
j\geq 1
\end{multline*}
(with the notational convention, followed also below,  that $ v^0_e
(x,\o)/\sqrt{\o(0,e)} :=0$ if $\o(0,e)=0$).
 Hence for Lebesgue a.a. $x\in \bbR^d$  the form
\begin{equation}\label{pietro}
 (\o ,e)\rightarrow  \frac{ v^0_e
(x,\o)}{\sqrt{\o(0,e)}}   -  \partial _e u^0 (x)
\end{equation}
belongs to $L^2_{\text{pot}}$.  Now fix $\varphi\in C_c ^\infty$ and
a function $\psi \in C(\O) $. Note that, if $\t_{x/\e} \oo (0,e)>0$,
then
$$
\e \nabla^\e_e \psi (\t_{x/\e}\oo)=\bigl[\nabla^{(\o)}_e
\psi\bigr](\t_{x/\e} \oo ) \,, \qquad x \in \e \cC (\oo)\,.
$$
Therefore we can write
$$
\e \nabla ^\e _e \left( \varphi(x) \psi (\t_{x/\e}\oo)\right) = \e
\left( \nabla ^\e _e  \varphi(x)\right) \psi ( \t_{x/\e+e}\oo)+
\varphi (x) \bigl[   \nabla^{(\o)}_e \psi\bigr]( \t_{x/\e}\oo) \,.
$$ Due to the above identity and (\ref{rinog}), taking the inner
product of (\ref{mammina}) with $\e \varphi (x) \psi (\t_{x/\e} \oo)
$, we obtain
\begin{equation*}
\begin{split}
& \e\l \int_{\bbR^d} u^\e_\oo (x) \varphi(x) \psi (\t_{x/\e}\oo)
\mu^\e_\oo(dx) +\\&  \frac{1}{2}\sum _{e\in \cB} \e \int_{\bbR^d}
\bigl(\t _{x/\e } \oo\bigr) (0,e)
 \nabla ^\e _e  u^\e_\oo (x) \,\nabla ^\e _e  \varphi(x)\, \psi (\t_{x/\e+e}\oo) \mu^\e_\oo(dx) +\\
& \frac{1}{2} \sum_{e\in \cB} \int_{\bbR^d} \bigl(\t _{x/\e }
\oo\bigr)  (0,e) \nabla ^\e _e  u^\e_\oo (x) \, \varphi(x)\,
\bigl[ \nabla^{(\o)}_e \psi \bigr] (\t_{x/\e}\oo) \mu^\e_\oo(dx) =\\
& \e \int_{\bbR^d} f^\e _\oo (x)   \varphi (x) \psi (\t_{x/\e}
\oo)\mu^\e_\oo (dx)   \,.
\end{split}
\end{equation*}
 We note that due to (\ref{mazinga}) all terms but the third one in the l.h.s.
are negligible as $\e\downarrow 0$ along the subsequence satisfying
(\ref{mora}) and (\ref{mela}). Hence, by definition of weak
two--scale limit and the trivial identity (\ref{daidai}), we
conclude that
\begin{equation*}
\sum _{e\in \cB} \int_{\bbR^d}dx \int_\O \mmu (d\o) \sqrt{\o (0,e)}
v^0 _e (x,\o) \varphi(x) \nabla^{(\o)}_e \psi(\o) =0\,.
\end{equation*}
The above identity can be rewritten as
\begin{equation*}
\sum _{e\in \cB} \int_{\bbR^d}dx \, \varphi(x)  \int_\O \mmu (d\o)
\o (0,e) \frac{ v^0 _e (x,\o)}{\sqrt{\o(0,e)}  } \nabla^{(\o)}_e
\psi (\o) =0\,.
\end{equation*}
 Due to the arbitrariness  of the
test functions $\varphi $, we get that for Lebesgue a.a. $x\in
\bbR^d$ it holds
\begin{equation*}
\sum _{e\in \cB}   \int_\O \mmu (d\o) \o (0,e) \frac{ v^0 _e
(x,\o)}{\sqrt{\o(0,e)}  } \nabla^{(\o)}_e \psi (\o) =0\,.
\end{equation*}
By a separability argument, this implies that $\theta_x \in
L^2_{\text{sol}} (M)$ for Lebesgue a.a. $x\in \bbR^d$. Since we know
that the form (\ref{pietro}) belongs to $L^2_{\text{pot}} (M)$, this
concludes the proof of (\ref{seppioline}).

\end{proof}

\section{Proof of Theorem \ref{prop_hom}
   }\label{kukakuka}

We start with a technical result, which could be proven in much more
generality:

\begin{Le}\label{regolarissimo} Fix $\o \in \O_*$. Let $h \in C  (\bbR^d)$
satisfy
\begin{equation}
|h(x) |\leq \frac{c}{1+|x|^{d+1}} \,, \qquad \forall x \in \bbR ^d
\,,
\end{equation}
 and suppose that \begin{equation}\label{mozart} L^2
(\mu^\e_\o) \ni h^\e_\o\rightarrow h \in L^2 (m dx )\,.
\end{equation}
Then
\begin{equation}\label{schubert}
\lim_{\e\downarrow 0 } \int_{\bbR^d} \bigl| h^\e_\o (x) - h (x)
\bigr|^2 \mu^\e_\o (dx ) =0\,.
\end{equation}
\end{Le}

\begin{proof}

Trivially, it is enough to prove the following limits
\begin{align}
& \lim_{\e\downarrow 0 } \int _{\bbR^d} h  (x) ^2 \mu^\e _\o (dx)
= m \int h(x)^2 dx  \,,\label{proko} \\
& \lim_{\e\downarrow 0 } \int _{\bbR^d}h^\e _\o (x)  h  (x) \mu^\e
_\o (dx) = m \int h(x)^2 dx\,,\label{proko1}
\\
&\lim_{\e\downarrow 0 } \int _{\bbR^d} h^\e_\o (x) ^2 \mu^\e _\o
(dx) = m \int h  (x)^2 dx\,.\label{proko2}
\end{align}
Since  $h \in L^2 (\mu^\e_\o)$,  the integrals in the l.h.s. of
(\ref{proko}) and (\ref{proko1})   are meaningful. Moreover, observe
that for each $\ell
>0$ one can find a function $g _\ell \in C  _c (\bbR^d)$ such
that $ h (x) = g_\ell  (x)  $ for any $x \in \bbR^d$ with $|x|\leq
\ell$, and $|g_\ell (x)|  \leq c / (1+|x|^{d+1})$.

\smallskip

In order to prove (\ref{proko}) we observe that
\begin{align}
& \left|
 \int _{\bbR^d} h (x) ^2 \mu^\e _\o (dx)
- \int_{\bbR^d} g_\ell (x)^2 \mu^\e_\o (dx) \right|\leq 2 \e ^d \sum
_{x \in \e \bbZ^d\,:\, |x|> \ell  }\frac{c}{1+|x|^{d+1}} \leq
c(\ell )\,, \label{cairo1}\\
& \left| \int _{\bbR^d} h (x) ^2 m  dx  - \int_{\bbR^d} g_\ell (x)^2
m dx \right|\leq 2 \int _{\{x \in \bbR^d\,:\, |x| >\ell \}}
 \frac{c}{1+|x|^{d+1}} m dx  \leq c(\ell ) \,,\label{cairo2}
\end{align}
for a positive constant $c(\ell)$ going to $0$ as $\ell \uparrow
\infty$. The above estimates (\ref{cairo1}) and (\ref{cairo2}), and
the limit
\begin{equation*}
\lim_{\e\downarrow 0 } \int _{\bbR^d} g _\ell (x) ^2 \mu^\e _\o (dx)
= m \int g_\ell  (x)^2 dx
\end{equation*}
(due to the definition of $\O_*$) allow to derive (\ref{proko}) by
taking the limit $\ell\uparrow \infty$.

\medskip

In order to prove (\ref{proko1}) we observe that
\begin{multline}\label{gianna1}
\left| \int _{\bbR^d}h^\e _\o (x)  h  (x) \mu^\e _\o (dx)  - \int
_{\bbR^d}h^\e _\o (x)  g_\ell  (x) \mu^\e _\o (dx)\right|\leq \\
 \|
h^\e _\o \|_{\mu^\e_\o} \| h - g_\ell  \|_{\mu^\e_\o}\leq c(\o )
 \left(2 \e ^d \sum
_{x \in \e \bbZ^d\,:\, |x|> \ell
}\frac{c^2}{\bigl(1+|x|^{d+1}\bigr)^2} \right)^{\frac{1}{2}}\leq
c(\o ) c(\ell)\,,
\end{multline}
and
\begin{multline}\label{gianna2}
\left| \int _{\bbR^d}  h^2  (x) m dx   - \int _{\bbR^d}h  (x) g_\ell
(x)  m dx \right|\leq  2 \|h\|_\infty  \int _{\{x \in \bbR^d\,:\,
|x|
>\ell \}}
 \frac{c}{1+|x|^{d+1}} m dx  \leq c(\ell )\,,
\end{multline}
 for a positive constant $c(\ell)$ going to $0$ as $\ell \uparrow
\infty$. Since $h^\e_\o \rightarrow h$ and $g_\ell \in C_c (\bbR^d)
$ we can conclude that
$$
\lim_{\e\downarrow 0 }\int _{\bbR^d}h^\e _\o (x)  g_\ell  (x) \mu^\e
_\o (dx) = m \int_{\bbR^d} h(x) g_\ell (x) dx \,.
$$
The above limit together with (\ref{gianna1}) and (\ref{gianna2})
implies (\ref{proko1}).

\medskip

Finally we observe that (\ref{proko2}) follows
  by applying (\ref{rimmelbis}) in
the definition of strong convergence with test functions $\varphi^\e
:= h^\e _\o$, $\varphi:= h$.

\end{proof}

We have now all the main tools in order to prove  Theorem
\ref{prop_hom}. We take $\o \in \O_*$, define $u^0, v^0_e$ as in
Lemma \ref{bracciano} and assume that $f^\e _\o \rightharpoonup
f_\o$. We want to prove that $u^0$ solves equation (\ref{caldo}) and
that (\ref{tommy}) holds.

First we observe that the weak two--scale convergence (\ref{mora})
implies  the weak convergence
\begin{equation}\label{alleviare}
 L^2 (\mu^\e_\o ) \ni u ^\e_\o   \rightharpoonup   u^0 \in L^2 (
m  dx)\,
\end{equation}
as $\e \downarrow 0$ along the subsequence of  Lemma
\ref{bracciano}.
 Taking the inner product of (\ref{mammina}) with a
test function $\varphi \in C_c ^\infty (\bbR ^d )$ and using
(\ref{rinog}), we get the identity
\begin{multline}
\l \int _{\bbR ^d} u^\e _\o (x) \varphi (x) \mu ^\e _\o (dx ) +
\frac{1}{2} \sum _{e \in \cB } \int _{\bbR^d}  \t_{x/\e} \o   (0,e)
\nabla ^\e_e u (x) \nabla ^\e_e \varphi  (x) \mu ^\e _\o (dx)=\\\int
_{\bbR ^d}  f^\e _\o (x) \varphi (x) \mu ^\e _\o (dx)\,.
\end{multline}
By taking the limit $\e \downarrow 0$ (along the subsequence of
Lemma \ref{bracciano}) and then dividing by $m$,  from the trivial
identity (\ref{daidai}),  the limit (\ref{mela}) in Lemma
\ref{bracciano}, the limit (\ref{alleviare}) and the hypothesis
$L^2(\mu^\e _\o)\ni f^\e_\o
 \strano    f \in L^2 ( m  dx)$
  we get
\begin{multline}\label{cristalli}
\l   \int _{\bbR ^d} u^0  (x) \varphi (x) dx + \frac{1}{2 m }\sum
_{e \in \cB } \int _{\bbR^d} dx\, \partial _e \varphi (x) \int_{\O}
\mu (d\o') \sqrt{ \o' (0,e) }v^0 _e (x,\o')=\\   \int _{\bbR ^d}
f(x) \varphi (x) dx \,.
\end{multline}
The second member in (\ref{cristalli}) can be rewritten as
\begin{equation}\label{turchia}
\frac{1}{2 m }\sum _{e \in \cB_* } \int _{\bbR^d} dx\, \partial _e
\varphi (x) \int_{\O} \mu (d\o')  \left[\o' (0,e) \theta_x (\o',e)
-\o' (0,-e)\theta_x (\o',-e)  \right]\,.
\end{equation}

Due to Lemma \ref{bracciano}, $u^0 \in H^1 (\bbR^d, dx)$.  Given $x
\in \bbR^d$ we consider the gradients
$$
\z(x)  := \nabla \varphi (x) =\bigl(
\partial_e \varphi  (x) \bigr) _{e\in \cB_*} \, , \qquad  \xi (x) : = \nabla  u^0 (x)
=\bigl(
\partial_e u^0 (x) \bigr) _{e\in \cB_*}\,
$$
(the definition is well posed for Lebesgue a.a.  $x\in \bbR^d$,
since $u^0 \in H^1 (\bbR^d, dx)$). Due to  Lemma \ref{bracciano}, we
know that
 for Lebesgue a.a. $x\in \bbR ^d$ the form $\theta_x$ defined in
 (\ref{mamma1}) coincides with
 the form $\p  w^{\xi(x)} $ (recall
definition (\ref{elefante})). Therefore, due to (\ref{zenfira}),
(\ref{totti}) and  (\ref{sensi}),   we can rewrite (\ref{turchia})
as
\begin{equation}\label{ankara}
\frac{1}{2m} (w^{\z(x)}, \p w ^{\xi(x)} ) _{L^2 (M)} = \frac{1}{2m}
(\z(x) , D \xi(x) ) =   ( \z(x), \cD \xi(x)  ) =   \bigl(\nabla
\varphi (x) ,\cD  \nabla u^0 (x)\bigr)\,.
\end{equation}
In conclusion, (\ref{cristalli}) reads
\begin{equation}\label{cristallibis}
\l   \int _{\bbR ^d} u^0  (x) \varphi (x) dx +   \int_{\bbR^d}
\bigl( \nabla u^0 (x), \cD \nabla  \varphi (x) \bigr) dx = \int
_{\bbR ^d} f(x) \varphi (x) dx \,.
\end{equation}
Hence, the function $  u ^0 $ of Lemma \ref{bracciano} is the
solution of equation (\ref{caldo}), which is unique (in particular
$u^0$ does not depend from $\o \in \O_*$).
Due to  Lemma \ref{lemmino} it is simple to verify that for each
sequence $\e_k \downarrow 0$ one can extract a sub--subsequence
$\e_{k_n}$ satisfying Lemma \ref{bracciano}.  Hence, by the previous
results, we conclude that  for each sequence $\e_k \downarrow 0$ one
can extract a sub--subsequence $\e_{k_n}$ such that
$$
 L^2 (\mu^{\e_{k_n} }_\o ) \ni u ^{\e_{k_n}} _\o   \rightharpoonup   u^0 \in L^2 (
 m  dx)\,,
$$
thus implying  that  the functions $u ^\e _\o\in L^2 (\mu ^\e _\o )
$ weakly converge to $ u^0 \in L^2 (m dx) $. This concludes the
proof of point (i).

\smallskip

In order to prove   the strong convergence of $u^\e_\o\in L^2
(\mu^\e_\o) $ to $u^0\in L^2 (mdx)$ in point (ii)  one can proceed
as in \cite{ZP}[Proof of Theorem 6.1]. We give the proof  for the
reader's convenience. Due to Lemma \ref{caldissimo} we only need to
prove that
\begin{equation}\label{sabato}
\lim_{\e\downarrow 0 } \int _{\bbR ^d} u^\e_\o  (x)^2 \mu^\e _\o
(dx) = m \int _{\bbR^d} u^0(x)^2  dx\,.
\end{equation}
To this aim, we  define $v^\e_\o$ as the solution in
$L^2(\mu^\e_\o)$ of the equation
\begin{equation}\label{patata1}
\l v^\e _\o - \cL ^\e _\o v^\e _\o = u^\e _\o\,.
\end{equation}
As already proven, $L^2(\mu^\e_\o) \ni u^\e _\o\strano u^0 \in L^2
(mdx) $. Hence, by applying point (i) of Theorem \ref{prop_hom}, we
can conclude that
\begin{equation}\label{cullando}
L^2(\mu^\e_\o) \ni v^\e _\o\strano v  \in L^2 (mdx) \,,
\end{equation}
where $v \in L^2(m dx )$ solves the equation
\begin{equation}\label{patata2}
\l v - \nabla \cdot ( \cD \nabla    v) = u^0\,.
\end{equation}
By taking the inner product of (\ref{mammina}) with $v^\e_\o$ and
then subtracting the identity obtained by taking the inner product
of (\ref{patata1}) with $u^\e_\o$, one obtains that
\begin{equation}\label{famelico}
(v^\e_\o , f^\e _\o )_{\mu^\e_\o} =  (u^\e_\o , u^\e _\o
)_{\mu^\e_\o}\,.
\end{equation}
Similarly, by taking the inner product of (\ref{caldo}) with $v $
and then subtracting the identity obtained by taking the inner
product of (\ref{patata2}) with $u^0$ one obtains that
\begin{equation}\label{famelicobis}
\int_{\bbR^d} v(x) f(x) dx = \int_{\bbR^d}  u^0(x)^2 dx \,.
\end{equation}
Since by assumption  $ L^2(\mu^\e_\o) \ni f^\e _\o\rightarrow  f \in
L^2 (mdx)$, from (\ref{cullando}) and the definition of strong
convergence we derive that
\begin{equation}
\lim_{\e\downarrow 0 } \int _{\bbR ^d} v^\e _\o (x) f^\e_\o (x)
\mu^\e _\o (dx) = m \int _{\bbR^d} v(x) f(x) dx \,.
\end{equation}
Due to (\ref{famelico}) and (\ref{famelicobis}), the above limit is
equivalent to (\ref{sabato}). As already mentioned, this limit and
Lemma \ref{caldissimo} imply (\ref{tommy2}).

\medskip

We finally prove point (iii). Let $f \in C  _c (\bbR^d)$ and define
$f^\e _\o $ as the function $f$ restricted on $\e \cC (\o)$. Then,
due to (\ref{presilla}), $ L^2 (\mu^\e _\o ) \ni f ^\e _\o
\rightarrow f\in L^2 (m dx )$. Due to point (ii) proven above, we
know that $ L^2 (\mu^\e _\o ) \ni u ^\e _\o \rightarrow u ^0 \in L^2
(m dx )$. Since the function $u^0$ solves (\ref{caldo}) with $f\in
C_c (\bbR^d) $, $u^0$ is continuous and decays fast at infinity (see
Exercise 3.13 in Section III.3 of \cite{RW}). In order to conclude
it is enough to apply Lemma \ref{regolarissimo}.

\section{Proof of Corollary  \ref{tg5}  }\label{kukakukabis}

 Due to a
generalization of the Trotter--Kato Theorem  \cite{ZP}[Theorem 9.2],
\cite{P}[Theorem 1.4], Theorem \ref{prop_hom} (ii) implies for each
$\o \in \O_*$  that
\begin{equation}\label{carlabruni} L^2
(\mu^\e _\o ) \ni P_{t,\o}^\e f^\e _\o \rightarrow P_t f \in L^2 (m
dx) \,,  \end{equation}
 whenever $L^2
(\mu^\e_\o) \ni f^\e_\o \rightarrow f\in L^2 (m dx )$. Since, it
holds $L^2(\mu^\e_\o ) \ni f \rightarrow f \in L^2 (m dx )$  for
each $f \in C_c (\bbR^d)$ and each $\o \in \O_*$, (\ref{carlabruni})
is verified by setting $f^\e_\o :=f$. This fact and  Lemma
\ref{regolarissimo} allow to derive (\ref{kerbala1}).

In order to conclude  we only need to derive (\ref{kerbala2}) from
(\ref{kerbala1}). To this aim, let $\L_\ell := [-\ell, \ell ]^d$,
$\ell>0$. We claim that, for any $\o\in \O_*$,  given any $f\in C_c
(\bbR^d)$ it holds
\begin{equation}\label{consulta}
\lim_{\e \downarrow 0 } \int _{\L_\ell^c } P^\e _{t, \o } f(x) \mu
^\e _\o (dx) = \int _{\L_\ell ^c } P_t f (x) m dx \,.
\end{equation}
Without loss of generality we can assume that $f \geq 0$.

Since
$$P\bigl[  X(t\e^{-2}|x)=z \bigr] =P\bigl[ X(t\e^{-2}| z)=x \bigr]\, \qquad \forall t\geq 0, \; \forall x,z \in \cC
(\o)\,,
$$
we can write
\begin{multline}
\int _{\bbR ^d } P^\e _{t, \o } f(x) \mu ^\e _\o (dx) = \e ^d \sum
_{x\in
\cC (\o ) } \sum_{z\in \cC (\o) } f(\e z) P(  X(t\e^{-2}|x)=z ) =\\
 \e ^d \sum _{x\in \cC (\o ) } \sum_{z\in \cC (\o) } f(\e z) P(  X(t\e^{-2}|z)=x ) =
\e ^d  \sum_{z\in \cC (\o) } f(\e z)\rightarrow m \int_{\bbR ^d }
f(z) dz \,.
\end{multline}
The above limit and the identity $\int_{\bbR ^d } P_t f(z) dz=
\int_{\bbR ^d } f(z) dz$, following from  the symmetry of $P_t$,
implies (\ref{consulta})  with $\L_\ell^c $ replaced by $\bbR^d$.
Therefore, in order to prove (\ref{consulta})   it is enough to show
that
\begin{equation}\label{sonnone}
\lim_{\e \downarrow 0 } \int _{\L_\ell} P^\e _{t, \o } f(x) \mu
^\e_\o (dx)  = \int _{\L_\ell } P_t f (x) m dx\,.
\end{equation}
To this aim we apply Schwarz inequality and obtain the bounds
\begin{multline}\label{joda}
\left| \int _{\L_\ell} P^\e _{t, \o } f(x) \mu ^\e_\o (dx)- \int
_{\L_\ell} P_t f (x)    \mu ^\e_\o (dx)\right| \leq \int _{\L_\ell}
\bigl| P^\e _{t, \o } f(x) -P_t f (x) \bigr|  \mu ^\e_\o (dx)\\ \leq
\mu ^\e_\o (\L_\ell )^{1/2}\left( \int _{\L_\ell} \bigl| P^\e _{t,
\o } f(x) -P_t f (x) \bigr|^2  \mu ^\e_\o (dx)\right)^{1/2}\,.
\end{multline}
Since  by (\ref{carlo})    $\mu ^\e_\o (\L_\ell )\rightarrow m (2
\ell) ^d $ for each $\o \in \O_*$, the above upper bound and
(\ref{kerbala1}) imply that the first member in (\ref{joda}) goes to
$0$ as $\e \downarrow 0$ for each  $\o\in \O_*$. To conclude the
proof of (\ref{sonnone}) it is enough to observe that for each $\o
\in \O_*$ the integral $\int _{\L_\ell} P_t f (x) \mu ^\e_\o (dx)$
converges to $m \int _{\L_\ell} P_t f (x) dx $ since $P_t  f$ is   a
regular function fast decaying to infinity (the proof follows the
same arguments used in order to check (\ref{proko})). This concludes
the proof of (\ref{consulta}).

\smallskip

Let us come back to (\ref{kerbala2}). For each $\ell >0 $  we can
bound
\begin{equation}
\begin{split}
& \int _{\bbR^d} \bigl |P_{t,\o}^\e f (x) - P_t
f(x) \bigr| \mu^\e_\o (dx) \leq \\
& \int _{\L_\ell } \bigl |P_{t,\o}^\e f (x) - P_t f(x) \bigr|
\mu^\e_\o (dx)+\int _{\L_\ell^c } P_{t,\o}^\e f (x)  \mu^\e_\o
(dx)+\int _{\L_\ell^c  }   P_t f(x) \mu^\e_\o (dx)\leq \\
& \mu ^\e_\o (\L_\ell )^{1/2}\left( \int _{\L_\ell } \bigl
|P_{t,\o}^\e f (x) - P_t f(x) \bigr|^2  \mu^\e_\o
(dx)\right)^{1/2}+\int _{\L_\ell^c } P_{t,\o}^\e f (x) \mu^\e_\o
(dx)+\int _{\L_\ell ^c } P_t f(x) \mu^\e_\o (dx)\,.
\end{split}
\end{equation}
Due (\ref{kerbala1}) and (\ref{consulta}), by taking $\e \downarrow
0$ we get that for each $\o \in \O_*$
$$
\limsup _{\e\downarrow 0 }\int _{\bbR^d} \bigl |P_{t,\o}^\e f (x) -
P_t f(x) \bigr| \mu^\e_\o (dx) \leq 2m \int _{\L_\ell^c }   P_t f(x)
dx\,.
$$
By the arbitrariness of $\ell$ in the above estimate one derives
(\ref{kerbala2}).

\appendix

\section{Proof of Lemma \ref{silurino}}
Let us suppose that $\hat \o_c$, $c>0$, stochastically dominates a
supercritical Bernoulli bond percolation and prove that hypotheses
(H2) and (H3) are satisfied. We call $\bbP$ the law of $\hat \o _c$
on $\{0,1\}^{\bbE_d}$. Due to Strassen Theorem, there exists a
probability measure $\cP$ on the product space
$\bbX:=\{0,1\}^{\bbE_d}\times\{0,1\}^{\bbE_d}$ such that (i) $ \o_1
(b) \geq \o _2 (b)$ for each $b \in \bbE_d$, for $\cP$ almost all
$(\o_1,\o_2) \in \bbX$, (ii) the marginal law of $\o_1$ is $\bbP$
and (iii) the marginal law of $\o_2$ is a Bernoulli bond percolation
with supercritical parameter $p>p_c$. It is well known (see
\cite{G}) that  $\o_2$ has a.s. a unique infinite cluster whose
complement has only connected components of finite cardinality. This
implies the same property  for the random field $\o_1$, thus
assuring that  the random field $\o$ fulfills hypothesis (H2).

\smallskip

Let us now consider  hypothesis (H3). We recall the variational
characterization of $\cD$:
\begin{equation}\label{jc1}
(a, \cD a) =\frac{1}{m} \inf _{\psi \in B (\O) }\left\{
 \sum _{e\in \cB_*} \int _\O \o (0,e)  ( a_e+ \psi (\t_e\o)-\psi (\o) ) ^2 \bbI _{0,e\in \cC (\o) }
  \bbQ ( d\o)\right\} \,.
\end{equation}
Since $\o(0,e)\geq c\, \hat \o_c (0,e)$, we obtain that
\begin{equation}\label{jc2}
 (a, \cD a) \geq
\frac{c}{m} \inf _{\psi \in B (\O) }\left\{
 \sum _{e\in \cB_*} \int _\O \hat \o _c (0,e)  ( a_e+ \psi (\t_e\o)-\psi (\o) ) ^2 \bbI _{0,e\in \cC (\hat \o_c) }
  \bbQ ( d\o)\right\} \,,
\end{equation}
where now $\cC(\hat\o_c)$ denotes the unique infinite cluster of
$\hat\o_c$ (due to the previous observations, the definition is well
posed a.s.). Given $\psi \in B(\O)$ we can write $\psi = f+g$ where
$f= \bbE(\psi |\cF)$ and $g=\psi-\bbE(\psi|\cF)$, $\bbE$ being the
expectation w.r.t. $\bbQ$ and $\cF$ being the $\s$--algebra
generated by the random variables $\hat \o_c (b)$, $b \in \bbE ^d$.
Since $\bbE (g |\cF)=0$, it is simple to check that
\begin{multline}\label{jc3}
\int _\O \hat \o _c (0,e)  ( a_e+ \psi (\t_e \o)-\psi (\o) ) ^2 \bbI
_{0,e\in \cC (\hat \o_c) }
  \bbQ ( d\o)=\\
\int _\O \hat \o _c (0,e)  \bigl \{[a_e+ f (\t_e\o)-f (\o) ] ^2+
[g(\t_e \o)-g(\o) ] ^2 \bigr \} \bbI _{0,e\in \cC (\hat \o_c) }
  \bbQ ( d\o)\,.
\end{multline}
This shows that the infimum in the r.h.s. of (\ref{jc2}) is realized
by   $\cF$--measurable functions. Therefore,
 \begin{multline}\label{jc4}
 \text{r.h.s. of (\ref{jc2})}  = \frac{c}{m} \inf _{\substack{\psi
\in B(\O): \\\psi=\psi(\hat \o_c) } }\left\{
 \sum _{e\in \cB_*} \int _\O \hat \o _c (0,e)  ( a_e+ \psi (\t_e\o)-\psi (\o) ) ^2 \bbI _{0,e\in \cC (\hat \o_c) }
  \bbQ ( d\o)\right\} =\\
    \frac{c}{m} \inf _{\substack{\psi \in B(\bbX)\\ \psi=\psi(\o_1)
} } \left\{
 \sum _{e\in \cB_*} \int _\bbX   \o_1 (0,e)  ( a_e+ \psi (\t_e \o_1 )-\psi (\o_1 ) ) ^2 \bbI _{0,e\in \cC (\o_1) }
  \cP ( d\o_1, d\o_2)\right\} \geq \\
   \frac{c}{m} \inf _{\psi \in (\bbX) } \left\{
 \sum _{e\in \cB_*} \int _\bbX   \o_1 (0,e)  ( a_e+ \psi (\t_e(\o_1,\o_2) )-\psi (\o_1,\o_2) ) ^2 \bbI _{0,e\in \cC (\o_1) }
  \cP ( d\o_1, d\o_2)\right\} \,.
\end{multline}
Above  we have used that the law $\hat \o_c$ and the marginal law of
$\o_1$ coincide.  Moreover, we recall that $B(\cdot)$ denotes the
space of bounded Borel functions on the given topological space.

 A this point, since $\o_1 (0,e) \geq \o_2 (0,e)$ and $\bbI_{0,e \in
 \cC(\o_1) } \geq \bbI _{0,e\in \cC (\o_2)}$
$\cP$--a.s., we can obtain another lower bound by substituting in
the last expectation $\o_1 (0,e)$ and $\cC(\o_1)$ with $\o_2 (0,e)$
and $\cC (\o_2)$ respectively. By taking the conditional expectation
w.r.t. to the $\s$--algebra generated by $\o_2$ and using the same
arguments as above, we derive that the last expression in
(\ref{jc4}) is bounded from below by
\begin{multline}\label{solare}
   \frac{c}{m} \inf _{\substack{\psi \in B(\bbX):\\ \psi=\psi(\o_2) } }
\left\{
 \sum _{e\in \cB_*} \int _\bbX   \o_2 (0,e)  ( a_e+ \psi (\t_e \o_2 )-\psi (\o_2) ) ^2 \bbI _{0,e\in \cC (\o_2) }
  \cP ( d\o_1, d\o_2)\right\}=\\
   \frac{c}{m} \inf _{\psi \in B(\{0,1\}^{\bbE _d} )   }
\left\{
 \sum _{e\in \cB_*} \int _{ \{0,1\}^{\bbE_d} }  \o (0,e)  ( a_e+ \psi (\t_e \o )-\psi (\o) ) ^2 \bbI _{0,e\in \cC (\o) }
  \bbP_p ( d\o)\right\}
   \,,
\end{multline}
where $\bbP _p$ is the Bernoulli bond percolation with parameter
$p>p_c$. In order to prove that the diffusion matrix $\cD$ is
positive defined, we only need to show that the r.h.s. of
(\ref{solare}) is positive for $a\not =0$. We point out that, apart
multiplicative factors, the last infimum in (\ref{solare}) equals
$(a, D_p a)$, $D_p$ being the diffusion matrix of the simple random
walk on the supercritical infinite cluster. One only needs to prove
the positivity of $D_p$. This result has been proven in
\cite{DFGW}[pages 828--838] in any dimension for $p> 1/ 2$ (see in
particular Remark 4.16 in \cite{DFGW}[page 837]). There the authors
are able to bound from below the diffusion matrix by means of the
effective conductivities  of suitable resistor networks (this
reduction works without any restriction on $p$). The positivity of
the effective conductivities is then derived by applying percolation
results originally proven for  $p>1/2$. These results have been
improved (see\cite{GM}[page 454] and references therein) and the
improvement allows to extend
  the  positive bound on the effective  conductivities  to all $p>p_c$.
  Other derivations of the positivity of $D_p$ can
be found in \cite{SS}, \cite{BB} and  \cite{MP}.

\smallskip
If $\hat \o$ is a  Bernoulli bond percolation with parameter
$p>p_c$, then for each $c>0$ the random field $\hat \o_c$ is a
Bernoulli bond percolation with parameter $p(c)$ such that $\lim
_{c\downarrow 0} p(c)=p$. Hence, taking $c>0$ small enough, we
obtain that hypotheses (H2) and (H3) are satisfied.

\smallskip
The last statement regarding the cases of $\cD$ diagonal or multiple
of the identity can be proved by the same arguments used in the
proof of Theorem 4.6 (iii) in \cite{DFGW}.

\bigskip \bigskip

\noindent {\sl Acknowledgements}.  The author   thanks  A.L.\,
Piatnitski for useful discussions during her visit to the
  Centre de Math\'ematiques et d'Informatique (CMI), Universit\'e
     de Provence, which she thanks  for the kind hospitality and the
  financial support.   Moreover, she
    acknowledges the financial support of GREFI--MEFI and
thanks the anonymous referees for  useful suggestions.


\begin{thebibliography}{99}


\bibitem[A]{A}  G.\ Allaire. {\em Homogenization and two--scale
convergence}. SIAM J. Math. Anal., {\bf 23}, 1482--1518 (1992).



\bibitem[B]{B}
P.\ Billingsley. {\em Convergence of probability measures.}  Second
edition.  J. Wiley, New York (1999).

\bibitem[BB]{BB} N.\ Berger, M.\ Biskup.
\emph{Quenched invariance principle for simple random walk on percolation clusters}.
 Probab. Theory Related Fields  {\bf 137},  no. 1-2, 83--120 (2007).

\bibitem[BP]{BP}
M.\ Biskup, T.M.\ Prescott. \emph{ Functional CLT for random walk
among bounded random conductances}. Electronic Journal of
Probability {\bf 12}, 1323--1348 (2007).


\bibitem[DFGW]{DFGW} A.\, De Masi, P.\,  Ferrari, S.\,  Goldstein, W.D.\, Wick. \emph{An invariance principle for reversible Markov processes. Applications to random motions in random
environments}. J. Statis. Phys. {\bf{55}} (3/4), 787--855 (1985).

\bibitem[D]{D} R.\, Durrett. {\em Ten Lectures on Particle Systems}. In  {\em Lecture Notes in Mathematics}, {\bf
1608},  Springer, Berlin  (1995).



\bibitem[F]{F} A.\, Faggionato. \emph{Bulk diffusion of 1D exclusion process with bond
disorder.} Markov Processes and Related Fields {\bf 13}, 519--542
(2007).



\bibitem[FJL]{FJL} A.\ Faggionato, M.\ Jara, C.\ Landim.
\emph{Hydrodynamic limit
 of one dimensional subdiffusive exclusion processes
  with random conductances}. To appear in Probab. Theory and Related Fields.


\bibitem[FM]{FM}
A.\ Faggionato, F.\ Martinelli. \emph{Hydrodynamic limit of a
disordered lattice gas}.  Probab. Theory and Related Fields  {\bf
127} (3), 535--608 (2003).

\bibitem[G]{G}  G.\ Grimmett. {\em Percolation.} Second edition. Springer,   Berlin (1999).


\bibitem[GM]{GM} G.\ Grimmett, J.\ Marstrand,  {\em The supercritical phase of percolation is
well behaved}. Proc. Royal Society (London) Ser. A. 4306, 429-457
(1990)

\bibitem[J]{J} M. Jara. \emph{Hydrodynamic limit for the simple
exclusion process on non--homogeneous graphs}. IMPA preprint.

\bibitem[JL]{JL} M. Jara, C. Landim. {\em  Nonequilibrium central limit theorem
  for a tagged particle in symmetric simple exclusion}.
   Annales de l'institut Henri Poincar\'{e} (B) Probabilit\'{e}s et Statistiques, {\bf 42},
   567--577, (2006).



\bibitem[KL]{KL}
C.\ Kipnis, C.\ Landim. \emph{ Scaling limits of interacting
particle systems}. Springer, Berlin (1999).


\bibitem[Ko]{Ko} S.M.\, Kozlov. \emph{The method of averaging and
walks in inhomogeneous environments}. Russian Math. Surveys {\bf
40}, 73--145 (1985)

\bibitem[Ku]{Ku} R.\, K\"{u}nnemann. \emph{The diffusion limit for
reversible jump processes on $\bbZ^d$ with ergodic random bond
conductivities}. Comm. Math. Phys. {\bf 90}, 27--68 (1983).


\bibitem[L]{L} T.M.\, Liggett. \emph{Interacting particle systems}.
 Springer,  New York  (1985).

\bibitem[M]{M} P.\ Mathieu. \emph{Quenched invariance principles for random walks with random
conductances}.  J. Stat. Phys. {\bf 130}, 1025-–1046 (2008).



\bibitem[MP]{MP}  P.\, Mathieu, A.L.\, Piatnitski. \emph{Quenched
invariance principles for random walks on percolation clusters}.
  Proceedings of the Royal Society A.   {\bf 463},   2287–-2307   (2007).

\bibitem[N]{N} K. Nagy. \emph{Symmetric random walk in random
environment}. Period. Math. Hung. {\bf 45}, 101--120 (2002).

\bibitem[Nu]{Nu} G.\ Nguetseng. \emph{A general convergence result
for a functional related to the thoery of homogenization}. SIAM J.
Math. Anal., {\bf 20}, 608--623 (1989).

\bibitem[P]{P} S.E.\ Pastukhova. \emph{On the convergence of
hyperbolic semigroups in a variable Hilbert space}.  J. Math. Sci.
(N.Y.), {\bf 127}, no. 5, 2263--2283 (2005).


\bibitem[PR]{PR} A.\, Piatnitski, E.\, Remy. \emph{Homogenization of
elliptic difference operators}. SIAM J. Math. Anal. {\bf 33}, 53--83
(2001)


\bibitem[Q1]{Q1}
J.\ Quastel. \emph{Diffusion in disordered media}. In
\emph{Proceedings in Nonlinear Stochastic PDEs} (T. Funaki and W.
Woyczinky, eds), Springer, New York, 65--79 (1996).


\bibitem[Q2]{Q2}
J.\ Quastel, \emph{Bulk diffusion in a system with site disorder.}
Ann. Probab. {\bf 34} (5),   1990--2036  (2006)



\bibitem[RW]{RW} L.C.G.\ Rogers, D.\ Williams. {\em Diffusions,
Markov processes and martingales}. Volume 1, Second edition,
Cambridge University Press, Cambridge (2005).



\bibitem[SS]{SS}  V.\ Sidoravicius, A.-S.\ Sznitman. {\em
Quenched invariance principles for walks on clusters of percolation
or among random conductances}. Probab. Theory Related Fields {\bf
129},  219–244 (2004).



\bibitem[Z]{Z}  V.V.\ Zhikov. \emph{On an extension of the method
of two--scale convergence and its applications}. (Russian)  Mat. Sb.
{\bf 191},  no. 7, 31--72 (2000);  translation in  Sb. Math. {\bf
191}, no. 7-8, 973--1014 (2000).


\bibitem[ZP]{ZP}  V.V.\, Zhikov, A.L.\, Pyatnitskii. \emph{Homogenization
of random singular structures and random measures}. (Russian) Izv.
Ross. Akad. Nauk Ser. Mat. {\bf 70}, no. 1, 23--74 (2006);
translation in Izv. Math. {\bf 70}, no. 1, 19--67 (2006).




\end{thebibliography}
\end{document}